\documentclass[aos,MSNbibl,nameyear,dvips]{arximspdf}
\doi{10.1214/12-AOS986} 
\volume{40}
\issue{2}
\pubyear{2012}
\firstpage{996}
\lastpage{1023}

\makeatletter
\newtheorem{theorem}{Theorem}
\newtheorem{lemma}[theorem]{Lemma}
\newproclaim{example}{Example}
\newtheorem{corollary}[theorem]{Corollary}

\newproclaim{remark}[theorem]{Remark}
\newproclaim{definition}[theorem]{Definition}
\newcommand{\m}[1]{\mathrm{#1}}
\newcommand{\bfn}{\mathbf{n}}
\newcommand{\bfm}{\mathbf{m}}
\newcommand{\bft}{\mathbf{t}}
\newcommand{\bfx}{\mathbf{x}}
\newcommand{\bfz}{\mathbf{z}}
\newcommand{\bfmu}{\bolds{\mu}}
\newcommand{\bftheta}{\bolds{\theta}}
\newcommand{\bfeta}{\bolds{\eta}}
\newcommand{\bflambda}{\bolds{\lambda}}
\newcommand{\bfgamma}{\bolds{\gamma}}
\newcommand{\eqref}[1]{(\ref{#1})}
\makeatother

\begin{document}
\begin{frontmatter}

\title{Maximum likelihood estimation in log-linear models\thanksref{T1}}
\runtitle{Maximum likelihood estimation in log-linear models}

\thankstext{T1}{Supported in part by NSF Grant DMS-06-31589.}

\begin{aug}
\author[A]{\fnms{Stephen E.} \snm{Fienberg}\ead[label=e1]{fienberg@stat.cmu.edu}}
\and
\author[B]{\fnms{Alessandro} \snm{Rinaldo}\corref{}\ead[label=e2]{arinaldo@cmu.edu}}

\runauthor{S. E. Fienberg and A. Rinaldo}
\affiliation{Carnegie Mellon University}
\address[A]{Department of Statistics\\
Machine Learning Department\\
Cylab\\
Heinz College\\
Carnegie Mellon University\\
5000 Forbes Avenue\\
Pittsburgh, Pennsylvania 15213\\
USA\\
\printead{e1}} 

\address[B]{Department of Statistics\\
Carnegie Mellon University\\
5000 Forbes Avenue\\
Pittsburgh, Pennsylvania 15213\\
USA\\
\printead{e2}}
\end{aug}

\received{\smonth{9} \syear{2011}}
\revised{\smonth{2} \syear{2012}}

%
\begin{abstract}
We study maximum likelihood estimation in log-linear models under
conditional Poisson sampling schemes. We derive necessary and
sufficient conditions for existence of the maximum likelihood estimator
(MLE) of the model parameters and investigate estimability of the
natural and mean-value parameters under a nonexistent MLE. Our
conditions focus on the role of sampling zeros in the observed table.
We situate our results within the framework of extended exponential
families, and we exploit the geometric properties of log-linear models.
We propose algorithms for extended maximum likelihood estimation that
improve and correct the existing algorithms for log-linear model
analysis.
\end{abstract}

%
\begin{keyword}[class=AMS]
\kwd[Primary ]{62H17}
\kwd[; secondary ]{62F99}.
\end{keyword}
\begin{keyword}
\kwd{Extended exponential families}
\kwd{extended maximum likelihood estimators}
\kwd{Newton--Raphson algorithm}
\kwd{log-linear models}
\kwd{sampling zeros}.
\end{keyword}

\end{frontmatter}

\section{Introduction}

Log-linear models are arguably the most popular and important
statistical models for the analysis of categorical data; see, for
example, \citet{BFH75}, \citet{CHRISTENSEN97}. These powerful models,
which include as special cases graphical models [see, e.g.,
\citet{LAU96}] as well as many logit models [see, e.g., \citet{AGR02},
\citet{BFH75}], have applications in many scientific areas, ranging
from social and biological sciences, to privacy and disclosure
limitation problems, medicine, data-mining, language processing and
genetics. Their popularity has greatly increased in the last decades
due to growing demands for analyzing databases taking the form of large
and sparse contingency tables, where most of the cell entries are very
small or zero counts. Despite the widespread usage of these models, the
applicability and statistical properties of log-linear models under
sparse settings are still very poorly understood. As a result, even
though high-dimensional sparse contingency tables constitute a type of
data that is common in practice (e.g.,\vadjust{\goodbreak} in sample survey applications),
their analysis remains exceptionally difficult; see \citet{EFJ07} for
such an example.


In this article we are concerned with statistical inference in
log-linear models of arbitrary dimension,
and, in particular, with conditions for the existence of the maximum
likelihood estimator, or MLE, of the
model parameters. In log-linear model analysis, virtually all
methodologies for assessment of fit, model
selection and interpretation are applicable and have theoretical
validity only provided that the MLE exists.
Though this may appear to be only a computational issue, in fact, when
MLE is not defined, the applicability of statistical procedures
routinely used by practitioners may no longer have a theoretical
justification and, at the very least, require alteration. The
statistical implications of a nonexistent MLE, some of which are
detailed below, are numerous and severe.
\begin{itemize}
\item Existence of the MLE is required to justify the use of large
sample $\chi^2$ approximations to numerous measures of goodness-of-fit
commonly utilized for model assessment and model selection; see, for
example, \citet{BFH75}, \citet{AGR02}, \citet{CR88}. When the MLE does
not exist, the standard regularity conditions used to derive such
approximations no longer hold. As we show below, under a~nonexistent
MLE, the model is not identifiable, the asymptotic standard
errors are not well defined and the number of degrees of freedom
becomes meaningless. Though existence of the MLE is by no means enough
to warrant the use of $\chi^2$ approximations, nonexistence will surely
make them inadequate.

%
\item Existence of the MLE is also needed to derive a limiting
distribution for the double-asymptotic approximations of the likelihood
ratio and Pearson's $\chi^2$ statistic for tables in which both the
sample size and the number of cells are allowed to grow unbounded, a
setting studied, among others, by \citet{MORRIS75}, \citet{HAB77} and
\citet{KOE86}; see also \citet{CR88}.
%
\item The issue of nonexistence is also important for Bayesian analysis
of log-linear models; see, for example, \citet{KB01}, \citet{MLD09},
\citet{DM10} and references therein. Indeed, we will demonstrate that
nonexistence of the MLE is due to the data not being fully informative
about the model parameters, and results in nonestimability of those
parameters. Since the nonexistence of MLEs
is due to insufficient data, it cannot be remediated. 
In particular, the use of Bayesian methods in cases in which the MLE is
nonexistent is equivalent to replacing the information content lacking
in the data with the information contained in the prior. 
Since for some parameters no learning from the data takes place, the
posterior distribution must be interpreted accordingly. Furthermore,
when one uses improper priors for the log-linear parameters, the
posterior may be also be improper when the MLE does not exist; see
\citet{JF}.
\end{itemize}

It has long been known [see, in particular, \citet{BIRCH63},
\citet{HAB74}, \citet{BFH75}] that the nonexistence of the MLE is
caused by sampling zeros. When certain patterns of zero counts occur in
the observed table, the log-likelihood function cannot be maximized by
any vector of finite norm. While for hierarchical log-linear models,
patterns of sampling zeros leading to null margins are well known to
cause nonexistence of the MLE, very little has been known or observed
about general patterns of sampling zeros associated with nonexistent
MLEs. The very few know examples described in \citet{HAB74},
\citet{ROY05} and \citet{DFRSZ09} suggest that nonexistence of the MLE
may occur in small tables, but is very likely to arise when the table
is large and sparse.

\citet{HAB74} first obtained necessary and sufficient conditions for
the existence of the MLE for log-linear models. \citet{MLE06} gave a
direct geometric interpretation of Haberman's conditons and proposed a~polynomial time algorithm for checking for the existence of the MLE.
\citet{AIK79} and \citet{VER92} refined Haberman's conditions by
recasting the problem within the frameworks of exponential families and
of generalized linear models, respectively. In fact, the issue of
nonexistence of the MLE is best dealt with using the general theory of
exponential families and, in particular, of extended exponential
families, originally put forward by \citet{BARN78} and then
\citet{BRW86}. See also the important work by \citet{CENCOV82}.
In a recent series of papers, \citeauthor{CM01} (\citeyear{CM01,CM03,CM05,CM08}) broadened significantly the notions of extended
exponential families and extended maximum likelihood estimation to
include very general settings under minimal assumptions. See, in
particular, Remark~5.9 in \citet{CM08}, which briefly point to the
connections with the theory of log-linear models. \citet{ERGM09} and
\citet{GEYER09} contain more specialized results directly relevant to
the log-linear settings. Adopting a different approach, \citet{LAU96}
{\it defined} the parameter space for log-linear models as the
point-wise limit closure of the log-linear model parameters, which he
calls the extended log-linear model, and effectively treats the MLE and
extended MLE as one entity. While this is theoretically convenient, the
issue of nonestimability of the model parameters is not resolved, and
the computation of the extended MLE is just as problematic. Finally,
\citet{NR11} provided asymptotic conditions under which, for a
hierarchical log-linear model, a penalized maximum likelihood estimator
based on the group-lasso penalty will return the correct model, with
high probability.


Despite the breadth of the cited literature, two key issues concerning
maximum likelihood estimation in log-linear models remain. First, the
properties of extended exponential families have not yet been
specialized to the case of log-linear models. In particular, direct
application of this theory does not yield, in general, usable
conditions for the existence of the MLE, and the identification of the
nonestimable log-linear parameters or of the patterns of zeros leading
to a nonexistent MLE are still open problems. Secondly, existing
theoretical results have not been incorporated yet in any numerical
algorithm for checking for existence of the MLE and for identifying
nonestimable parameters. Consequently, virtually all statistical
software currently available to practitioners is flawed, to the point
that nonexistence of the MLE can be detected only by monitoring whether
the algorithm used to optimize the log-likelihood function fails to
converge, or converges slowly or becomes unstable; see, for example,
\citet{ROY05}. Consequently, results and decisions stemming from the
statistical analysis of contingency tables containing substantial
numbers of zero counts can be seriously compromised.

In this article we attempt to rectify these problems. 
Our contributions are two-fold:
\begin{itemize}
\item From a theoretical standpoint, we derive necessary and sufficient
conditions for existence of the MLE that are broadly applicable to a
variety of sampling schemes and amenable to computations. Ultimately,
these conditions amount to checking whether the observed sufficient
statistics lie on the boundary a polyhedral cone, called the {\it
marginal cone}; see \citet{MLE06}. When the MLE does not exist, we
specialize the theory of extended exponential families to characterize
the estimability of the natural and mean-value parameters of the
log-linear models. To this end, we focus on discrete exponential
families with polyhedral convex support [see \citet{ERGM09},
\citet{GEYER09}], and rely significantly on tools from polyhedral
geometry.
\item From a practical viewpoint, we develop algorithms for extended
maximum likelihood estimation that are applicable to large tables. Our
procedures will allow one to (i) detect nonexistence of the MLE and
(ii) identify and estimate all the parameters that are in fact
estimable. Overall, our algorithms correct and improve over many
existing software for log-linear model analysis. Due to space
constraints, a detailed description of these algorithms is contained in
the supplementary material [\citet{FRsupp}].
\end{itemize}

\subsection*{Notation}

We let $\mathcal{I}$ be a finite set of indices or cells, representing
the support of a discrete distribution, such as the joint distribution
of a set of categorical variables. We set $I = |\mathcal{I}|$, where
$|B|$ is the cardinality the set $B$. We denote by
$\mathbb{R}^{\mathcal{I}}$ be the vector space of real-valued functions
on $\mathcal{I}$, and $\mathbb{R}^{\mathcal{I}}_{\geq0}$ and
$\mathbb{N}^{\mathcal{I}}$ its subset of nonnegative functions and
nonnegative integer-valued functions, respectively. For vectors
$\mathbf{x}$ and $\mathbf{y}$, $(\mathbf{x}, \mathbf{y}) =
\mathbf{x}^\top\mathbf{y} $ represents their inner product and $\|
\mathbf{x}\| = \sqrt{\bfx^\top \bfx}$ the corresponding Euclidean norm.
If $\bfx\in \mathbb{R}^\mathcal{I}$, we denote by $\mathbf{x}(i)$ the
value corresponding to the $i$th coordinate of $\mathbf{x}$ and by
$\operatorname{supp}(\mathbf{x}) = \{ i \dvtx \mathbf{x}(i) \neq0\}$
the set of coordinates of $\mathbf{x}$ with nonzero values. We take
functions and relations on vectors component-wise, for example, for
$\mathbf{x} \in \mathbb{R}^{\mathcal{I}}$, $\exp(\mathbf{x}) = \{
e^{\mathbf{x}(i)} \dvtx  i \in\mathcal{I} \}$.

For a nonempty subset $\mathcal{F} \subseteq\mathcal{I}$, we let
$\pi_\mathcal{F} \dvtx  \mathbb{R}^{\mathcal{I}} \rightarrow
\mathbb{R}^{\mathcal{F}}$ the coordinate projection map given by $ \{
\mathbf{x}(i) \dvtx  i \in\mathcal{I} \} \mapsto\{ \mathbf{x}(i) \dvtx
i \in\mathcal{F} \} $ and, for any $S \subset\mathbb{R}^{\mathcal{I}}$,
we set $\pi_\mathcal{F}(S) = \{ \pi_{\mathcal{F}}(\mathbf{x}),
\mathbf{x} \in S\}$. If $\mathcal{M}$ is a linear subspace, we denote
by $\mathcal{M}^\bot$ its orthogonal complement and by
$\Pi_{\mathcal{M}}$ the orthogonal projector into $\mathcal{M}$. If
$\mathcal{N}$ is another linear subspace contained in $\mathcal{M}$, we
write $\mathcal{M} \ominus\mathcal{N}$ for the subspace $\mathcal{M}
\cap \mathcal{N}^\bot$.

For a matrix $\m{A}$, $\mathcal{R}(\m{A})$ denotes its column range and
$\operatorname{kernel}(\m{A})$ its null space. If the rows of $\m{A}$
are indexed by $\mathcal{I}$, and $\mathcal{F}$ is a nonempty subset of
$\mathcal{I}$, $\m{A}_{\mathcal{F}}$ is the submatrix of $\m{A}$
comprised of the rows with indexes in $\mathcal{F}$. We write
$\operatorname{cone}(\m{A})$ for the polyhedral cone spanned by columns
of $\m{A}$ and $\operatorname{conv}(\m{A})$ for the polytope consisting of the
convex combinations of its columns. Similarly, for a set $S$,
$\operatorname{conv}(S)$ is the convex hull of all its points. For a
polyhedron $P$, we write its relative interior as
$\operatorname{ri}(P)$.

\setcounter{footnote}{1}

\section{Log-linear models, sampling schemes and exponential
families}\label{sec:samplsch}

Log-linear model analysis is concerned with the study of discrete
probability distributions over a finite set $\mathcal{I}$, whose
elements will be referred to as {\it cells}. These distributions are
assumed to form an exponential family of probabilities $\{ P_{\bfeta},
\bfeta\in\mathbb{R}^d \}$ with densities with respect to the
counting measure on $\mathcal{I}$ of the form
\begin{equation}\label{eq:exp}
p_{\bfeta}(i) = P_{\bfeta}(\{ i \}) = \exp \{  ( \bfeta, \mathbf{a}_i
 ) - \phi(\bfeta)  \},\qquad \bfeta\in\mathbb{R}^d,
\end{equation}
where each $\mathbf{a}_i$ is a nonzero vector in $\mathbb{R}^d$, and
$\phi(\bfeta) = \log ( \sum_i \exp \{  (\bolds{\eta}, \mathbf{a}_i
 )  \}  )$ is the log-partition function. The $I \times
d$ matrix $\m{A}$, whose $i$th row is the vector~$\mathbf{a}^\top_i$,
is called the {\it design matrix}.\footnote{It is easy to see that
design matrices are not uniquely determined: if $\m{A}_1$ and $\m{A}_2$
are two matrices of dimensions $I \times d_1$ and $I \times d_2$,
respectively, and with identical row spans, then they parametrize the
same statistical model.}

Suppose we observe a sample of $N$ independent and identically
distributed realizations from an unknown distribution satisfying
\eqref{eq:exp}, where the data take the form of an unordered sequence
of random cells $(L_1, \ldots, L_N)$, with $L_j \in\mathcal{I}$ for
each $j$, and where $N$ too can be random. The observed cells are then
cross-classified into a random integer vector $\mathbf{n} \in
\mathbb{N}^{\mathcal{I}}$, called a {\it a contingency table}, with
$\mathbf{n}(i) = |\{ j \dvtx  L_j = i \}|$, for all $i \in\mathcal{I}$.


Traditionally, log-linear model analysis is not directly concerned with
the natural parameters $\bfeta$ in \eqref{eq:exp}, but rather with the
unknown expected value $\mathbf{m} : = \mathbb{E}[\mathbf{n}]$ of the
resulting contingency table, under the provision that \mbox{$\mathbf{m}(i) >
0$} for each $i$. In detail, letting $\mathcal{M} \subset
\mathbb{R}^{\mathcal{I}}$ be the linear subspace spanned by the rows of
the design matrix $\m{A}$, the ensuing log-linear model is predicated
on the condition that $\bolds{\mu} : = \log(\mathbf{m}) \in\mathcal
{M}$. In particular, log-linear models are typically defined as
statistical models for the distribution of the random table
$\mathbf{n}$ indexed by the points in the linear subspace
$\mathcal{M}$.



The distribution of the table $\mathbf{n}$ depends on the sampling
scheme used during the data collection process. In this article, we
study sampling schemes based on linear restrictions on $\mathbf{n}$,
known as conditional Poisson sampling schemes, introduced in
\citet{HAB74}, Chapter 1. Specifically, let $\mathcal{N}$ be a given
$m$-dimensional linear subspace of $\mathcal{M}$, which we will refer
to as the {\it sampling subspace}, and $\mathbf{c}$ a known vector in
$\mathbb{R}^{\mathcal{I}}$. The corresponding conditional sampling
\mbox{Poisson} scheme prescribes that the distribution of $\mathbf{n}$ is
given by the conditional distribution of $I$ independent Poisson random
variable $\{ \mathbf{n}(i), i \in\mathcal{I}\}$ with mean parameters
$\{ \mathbf{m}(i) = \exp (\bolds{\mu}(i)  ), i \in\mathcal{I}\}$, where
$\bolds{\mu} \in\mathcal{M}$, given that $\Pi_{\mathcal{N}} \mathbf{n}
= \mathbf{c}$.
This type of data sampling includes the most commonly used sampling
schemes, described below.
\begin{itemize}

\item\textit{Poisson sampling scheme.} The sampling subspace is
$\mathcal{N} = \{ \mathbf{0} \}$. Thus, there are no restrictions on
$\mathbf{n}$, which is a random vector comprised of independent Poisson
random variables with mean $\mathbf{m}$. The log-likelihood function is
given by
\begin{equation}\label{eq:pois.loglik.mu}
\ell^P(\bolds{\mu}) = (\mathbf{n},\bolds{\mu})
-\mathbf{1}^\top\exp(\bolds {\mu}) - \sum_{i} \log\mathbf{n}(i)!,\qquad
\bolds{\mu} \in\mathcal{M}.
\end{equation}

\item\textit{Product multinomial and multinomial sampling schemes.} Let
$\mathcal{B}_1, \ldots, \mathcal{B}_m$ be a partition of $\mathcal{I}$.
Under the product multinomial sampling, the conditional distribution of
the cell counts $\mathbf{n}$ is the product of $m$ independent
multinomials of sizes $N_j$, $j=1,\ldots,m$, each supported on the
corresponding class $\mathcal{B}_j$. Formally, let $\bolds{\chi}_j$ be
the indicator function of $\mathcal{B}_j$, where $\bolds{\chi}_j(i)$ is
$1$ if $i \in\mathcal{B}_j$ and $0$ otherwise,
%
and define $\mathcal{N}$ to be the $r$-dimensional subspace spanned by
the orthogonal vectors $(\bolds{\chi}_1,\ldots,\bolds{\chi}_r)$. The
product multinomial sampling constraints are of the form $(\mathbf{n},
\bolds{\chi}_j) = N_j$, for known integer constants $N_j$.
The log-likelihood function is [see \citet{HAB74}, equation~1.51]
\begin{eqnarray}\label{eq:multprodloglik}
&&\tilde{\ell}^M(\bolds{\mu}) = \sum_{j=1}^r  \biggl( \sum_{i \in
\mathcal{B}_j} \mathbf{n}(i)\log\frac{\mathbf{m}(i)}{(\mathbf{m},
\bolds{\chi }_j)} + \log N_j! - \sum_{i \in\mathcal{B}_j}
\log\mathbf{n}(i)!  \biggr),\nonumber
\\[-8pt]\\[-8pt]
&&\eqntext{\bolds{\mu} \in\mathcal{M},}
\end{eqnarray}
where $\mathbf{m} = \exp(\bolds{\mu})$. Because of the sampling
constraints, $\tilde{\ell}^M$ is well defined only on the subset of
$\mathcal{M}$,
\begin{equation}\label{eq:Mtilde}
\widetilde{\mathcal{M}} : = \{\bolds{\mu} \in\mathcal{M} \dvtx
(\bolds{\chi}_j, \exp(\bolds{\mu})) = N_j, j=1,\ldots,r\},
\end{equation}
which is  is neither a vector space nor a
convex set. We give a more convenient parametrization below in Lemma
\ref{lem:mpp}. The multinomial scheme is a~special case of product
multinomial schemes, corresponding to the trivial one-class partition
of $\mathcal{I}$ with indicator function $\mathbf{1}$. In this case,
$\mathbf{n}$ has a~multinomial distribution with size $N =
(\mathbf{1},\mathbf{n}) = (\mathbf{1},\bfm)$ and cell probabilities
$\mathbf{m}/N$.

\item\textit{Poisson--multinomial sampling schemes.} This sampling
scheme is a combination of the previous two schemes. For a given
partition $\mathcal{B}_1, \ldots, \mathcal{B}_m$ of $\mathcal{I}$, the
sampling constraints are of the form $(\mathbf{n}, \bolds{\chi}_j) =
N_j$ for $j=1,\ldots,m-1$, with the counts for the cells in the set
$\mathcal{B}_m$ left unconstrained; see \citeauthor{LANG04}
(\citeyear{LANG04,LANG05}).
\end{itemize}

As  is customary, we assume throughout that the sampling subspace
$\mathcal{N}$ is strictly contained in $\mathcal{M}$. The case
$\mathcal{N} = \mathcal{M}$ is practically uninteresting, as the
resulting sampling constraints would fix the value of the sufficient
statistics so that the conditional distribution of $\mathbf{n}$ will
not depend on the model parameters. We treat the case $\mathcal{N}
\not\subset\mathcal{M} $ in the supplementary material
[\citet{FRsupp}].

 We now derive the equivalent exponential family
representation for log-linear models under conditional Poisson schemes.
To this end, we will express the sampling constraints in a different,
but equivalent, form. Let $ (\mathbf{v}_1, \ldots, \mathbf{v}_m
 )$ be any set of $m$ vectors spanning $\mathcal{N}$ and such that
$(\mathbf{v}_j,\mathbf{c}) = 1$ for all $j$. Then, the sampling
constraints take the form
\[
\m{V}^\top\mathbf{n} = \mathbf{1},
\]
where $\m{V}$ is the $I \times m$ matrix whose $j$th column is
$\mathbf{v}_j$. Accordingly, we denote with
\[
S(\m{V}) : = \{\mathbf{x} \in\mathbb{N}^{\mathcal{I}}\dvtx  \m{V}^\top
\mathbf{x} = \mathbf{1} \}
\]
the set of all possible tables compatible with the sampling constraints
specified by~$\m{V}$. Let $\nu$ be the finite measure on
$\mathbb{N}^\mathcal{I}$ given by\footnote{This particular choice of
the dominating measure will lead to Poisson and product multinomial
likelihoods. More generally, much of our analysis carries over with
other choices of dominating measure, for example, the ones for which
conditions (A1)--(A4) in \citet{ERGM09} hold.}
\[
\nu(\bfx) : = \prod_{i \in\mathcal{I}}\frac{1}{\bfx(i)!},\qquad
\bfx\in\mathbb{N}^\mathcal{I}.
\]
For a conditional Poisson scheme defined by $\m{V}$, let $\nu_{\m{V}}$
be the restriction of~$\nu$ on $S(\m{V})$, that is, $\nu_{\m
{V}}(\bfx)
: = 1_{x \in S(\m{V})} \nu(\bfx),$ with $\bfx\in
\mathbb{N}^\mathcal{I}$.

It is easy to see that the conditional distribution of the table
$\bfn$, given the sampling constraints determined by $\m{V}$, is the
exponential family of distributions with base measure $\nu_{\m{V}}$,
sufficient statistics $\m{A}^\top\bfx$, natural parameter space
$\mathbb{R}^d$ and densities given by
\begin{equation}\label{eq:exp.A}
p_{\bftheta}(\bfx) = \exp \{ (\m{A}^\top\bfx,\bftheta) - \psi(\bftheta)
 \},\qquad \bfx\in S(\m{V}), \bftheta\in \mathbb{R}^d,
\end{equation}
where $\psi(\bftheta) = \log ( \int_{S(\m{V})} \exp \{ (\m{A}^\top
\bfx,\bftheta)  \} \,d \nu_{\m{V}}(\bfx)  )$. This exponential
representation is not the most parsimonious from the viewpoint of
sufficiency. Indeed, let $\mathcal{T} = \{ \bft\in \mathbb{R}^d \dvtx
\bft= \m{A}^\top\bfx, \bfx\in S(\m{V}) \}$ be the image of $\m{A}$ and
$\mu_{\m{V}} = \nu_{\m{V}}\m{A}^{-1}$ be the measure induced by
$\m{A}$. Then, by standard arguments [see, e.g., \citet{BRW86}], the
distributions of the sufficient statistics $\bft= \m{A}^\top\bfn$ also
form an exponential family, with density with respect to the base
measure $\mu_{\m{V}}$ given by
\begin{equation}\label{eq:qtheta}
q_{\bftheta}(\bft) = \exp \bigl( (\bft,\bftheta) - \psi(\bftheta )
\bigr),\qquad \bft\in\mathcal{T}, \bftheta\in\mathbb{R}^d,
\end{equation}
the same the log-partition function $\psi$ and natural parameter space
as in the original family.

It is now easy to see that the exponential family parametrization and
the log-linear parametrization are equivalent. Indeed, for any $\bfn$
and $\bft$ such that $\bft= \m{A}^\top\bfn$, and for any $\bftheta
\in\mathbb{R}^d$, the identity
\[
(\bft,\bftheta) = (\m{A}^\top\bfn,\bftheta) = (\bfn,\m{A}
\bftheta) = (\bfn,\bfmu),
\]
where $\bfmu= \m{A} \bftheta\in\mathcal{M}$, implies these models can
be equivalently parametrized by the linear subspace $\mathcal{M}$. If
$\m{A}$ is of full rank, then the map $\bftheta\mapsto\m{A} \bftheta$
is an isomorphism between $\mathbb{R}^d$ and $\mathcal{M}$, while if $d
> \operatorname{dim}(\mathcal{M})$, the natural parametrization is redundant
and, in fact, nonidentifiable.

Throughout this article, will impose the following assumptions. Let
$\m{V}$ be the matrix specifying the conditional Poisson sampling
scheme.

\begin{enumerate}[(A2)]
\item[(A0)] {\it Nontriviality}: the set $S(\m{V})$ is nonempty.
\item[(A1)] {\it Exhaustive sampling condition}: there does not exist
any vector $\bolds{\gamma} \in \mathcal{N}^\bot\setminus\{ 0 \}$, such
that $(\bolds{\gamma},\bfn )$ is constant almost everywhere with
respect to $\nu_{\m{V}}$. In particular, for no cell $i \in
\mathcal{I}$, $\mathbf{n}(i) = 0$, almost everywhere $\nu_{\m{V}}$.
\item[(A2)] {\it Integrality assumption}: $\{ \bfx\in\mathbb
{R}^{\mathcal{I}}_{\geq0} \dvtx \m{V} \bfx= \mathbf{1}\} =
\operatorname{conv} ( S(\m{V})  )$.
\end{enumerate}

Assumption (A1) guarantees that no linear constraints hold, other than
the ones specified by $\mathcal{N}$, and it prevents the sampling
constraints from introducing structural zeros. Even though we can
easily extend our analysis to deal with structural zeros, we do not
provide the details here. Assumption~(A2) is technical, and it is used
in Theorem \ref{thm:ex.mle} below to unify the conditions for existence
of the MLE across different sampling schemes. If (A2) is not in effect,
checking for existence of the MLE can become computationally
infeasible, depending on $\m{V}$. The Poisson, product multinomial and
Poisson--multinomial schemes automatically satisfy (A2).



\subsection{The effects of sampling constraints}

We conclude this section by studying the effect of the sampling
constraints on the estimability of the natural and log-linear
parameters. We show that imposing linear sampling restrictions results
in nonidentifiability of the corresponding natural exponential family
\eqref{eq:exp.A}, to the extent that only certain linear combinations
of the natural parameters, which depend only on the subspace
$\mathcal{N}$, are estimable. For the log-linear parameters, only
$\Pi_{\mathcal{M} \ominus\mathcal{N}} \bfmu$ is estimable, which
implies that the number of estimable parameters is
$\operatorname{dim}(\mathcal{M} \ominus\mathcal{N}) = d - m$.

We define the following equivalence relation on $\mathbb{R}^d$: for
$\bftheta_1,\bftheta_2 \in\mathbb{R}^d$, $\bftheta_1
\stackrel{\mathcal{N}}{\sim}\bftheta_2$ if and only if $\bftheta_1 -
\bftheta_2 \in\mathcal{Z}$, where
\begin{equation}\label{eq:Z}
\mathcal{Z} : = \{ \bolds{\zeta} \in\mathbb{R}^d \dvtx \m{A} \bolds
{\zeta} \in\mathcal{N}\}.
\end{equation}
%
For any $\bftheta\in\mathbb{R}^d$, we then write
$\bftheta_{\mathcal{N}} : = \{ \bftheta^* \dvtx \bftheta
\stackrel{\mathcal{N}}{\sim}\bftheta^* \}$ for the equivalence class
containing $\bftheta$, and $\Theta_{\mathcal{N}} : = \{
\bftheta_{\mathcal{N}} , \bftheta\in\mathbb{R}^d \}$ for the set of
equivalent classes corresponding to the equivalence relation
$\stackrel{\mathcal{N}}{\sim}$. For simplicity, below we assume that
the matrices $\m{A}$ and $\m{V}$ are of full rank, but the same
conclusions hold with $d$ replaced by $\operatorname{rank}(\m{A})$.


\begin{lemma}\label{lem:thetaN}
Consider the exponential family \eqref{eq:exp.A}, with $\m{A}$ of full
rank $d$, and suppose that conditions \textup{(A0)} and \textup{(A1)}
hold.
\begin{enumerate}[(ii)]
\item[(i)] 
The set $\Theta_{\mathcal{N}}$ is a vector space of dimension $d-m$
isomorphic to $\mathcal{M} \ominus\mathcal{N}$, and is comprised of
parallel $m$-dimensional affine subspaces of $\mathbb{R}^d$.
\item[(ii)] The family is nonidentifiable: any two points $\bftheta_1
\stackrel{\mathcal{N}}{\sim}\bftheta_2$ specify the same distribution.
In fact, this family is parametrized by $\Theta_{\mathcal{N}}$, or,
equivalently, by $\mathcal{M} \ominus\mathcal{N}$. Therefore, it is of
order $d - m$.
\end{enumerate}
\end{lemma}

Using standard minimality arguments, nonidentifiability of the natural
parameters can be easily resolved by redefining a smaller exponential
family of order $d - m$ using as a new design matrix any full-rank
matrix whose column span is $\mathcal{M} \ominus\mathcal{N}$; for this
fully-identifiable family, the natural parameter space is
$\mathbb{R}^{d-m}$. Concretely, we assume, without loss of generality,
that the matrix $\m{A}$ is of the form
\begin{equation}\label{eq:A}
\m{A} =  ( \m{B}\; \m{V} ),
\end{equation}
where $\m{V}$ is the $I \times m$ matrix of sampling restrictions
whose rows span $\mathcal{N}$ and $\m{B}$ is a $I \times(d - m)$
matrix whose row space is $\mathcal{M} \ominus\mathcal{N}$. Then,
replacing~$\m{A}$ with $\m{B}$ in \eqref{eq:exp.A} will produce a full
and minimal exponential family.

To illustrate this point, we show that the log-likelihood function
\eqref{eq:multprodloglik} for the product multinomial sampling scheme
can be more conveniently parametrized by $\mathcal{M} \ominus
\mathcal{N}$ instead of the nonconvex set $\widetilde{\mathcal{M}}$.
For any $\bolds{\beta} \in\mathcal{M} \ominus\mathcal{N}$, let
\begin{equation}\label{eq:multpar}
\ell^M(\bolds{\beta}) : = (\mathbf{n}, \bolds{\beta}) - \sum_{j=1}^m
N_j \log(\exp(\bolds{\beta}),\bolds{\chi}_j) - \sum_{i \in\mathcal{I} }
\log\bfn(i) !.
\end{equation}

\begin{lemma}\label{lem:mpp}
The sets $\mathcal{M} \ominus\mathcal{N}$ are $\widetilde{\mathcal{M}}$
homeomorphic and, for each pair of homeomorphic vectors
$\bfmu\in\widetilde{\mathcal{M}}$ and $\bolds{\beta} \in\mathcal{M}
\ominus\mathcal{N}$, $\tilde{\ell}_{\mathcal{L}}(\bolds{\mu}) =
\ell^M(\bolds{\beta}) $.\vadjust{\goodbreak}
\end{lemma}

The form of the likelihood in \eqref{eq:multpar} is better suited for
computations, as we show in \citet{FRsupp}.

Under the conditions of Lemma \ref{lem:thetaN}, the Fisher information
matrix at $\bftheta$ has rank $d-m$, for each $\bftheta\in
\mathbb{R}^d$. To see this, notice that the the Fisher information
matrix $I(\bftheta) $ at $\bftheta$ is
$\operatorname{Cov}_{\bftheta}(\m{A}^\top\bfn)$, where
$\operatorname{Cov}_{\bftheta}$ denotes the covariance operator
evaluated using the distribution parametrized by $\bftheta$. Then, for
any~$\bolds{\zeta}$ in the set $\mathcal{Z}$ defined in \eqref{eq:Z},
the linear form $(\m{A}^\top\bfn,\bolds{\zeta})$ is constant almost
everywhere and therefore has zero variance. This is equivalent to
$\bolds{\zeta}^\top I(\bftheta) \bolds{\zeta} = 0,$ so that
$\operatorname{rank}(I(\bftheta)) =
\operatorname{dim}(\mathcal{Z}^\bot) = d - m$, for all $\bftheta$.

\section{Theory of maximum likelihood estimation}\label{sec:mle}

We now provide a systematic treatment of maximum likelihood estimation
for the natural and log-linear parameters, within the framework of the
theory of discrete extended exponential families with linear sufficient
statistics. We refer the reader to \citet{BARN78} and \citet{BRW86} for
classic references and \citeauthor{CM01} (\citeyear{CM01,CM03,CM05,CM08}) for advanced treatments. In our setting, \citet{GEYER09}
and \citet{ERGM09} are particularly relevant. For the reader's
convenience, we briefly review the aspects of this theory that are
relevant to our problem in Appendix \ref{app:ext.expo}.


\subsection{Existence of the MLE}\label{sec:existence}

We prove a general necessary and sufficient condition for existence of
the MLE that applies to any conditional Poisson sampling scheme
satisfying assumptions (A0)--(A2). Unlike existing results, these
conditions directly translate into usable algorithms for checking for
the existence of the MLE, as described in \citet{FRsupp}.




For any design matrix $\m{A}$, we denote by $C_{\m{A}} : =
\operatorname{cone}(\m{A}^\top)$ the polyhedral cone spanned by the
rows of $\m{A}$. Following \citet{MLE06}, we call $C_{\m{A}}$ the {\it
marginal cone} of $\m{A}$.

\begin{theorem}\label{thm:ex.mle}
Assume conditions \textup{(A0)--(A2)} and let $\m{A}$ be any matrix
with column span $\mathcal{M}$. The MLE of $\bftheta_{\mathcal{N}}$
(or, equivalently, of $\Pi_{\mathcal{M} \ominus\mathcal{N}} \mu$)
exists and is unique if and only if $\bft= \m{A}^\top\bfn\in
\operatorname{ri}(C_\m{A})$.
\end{theorem}

This result is a nontrivial application of a well-known result about
existence of MLE in exponential families (viz., Theorem \ref{thm:mle}
in Appendix \ref{app:ext.expo}), and it subsumes previous results of
\citet{HAB74} and \citet{MLE06}, because it provides a unified
condition that applies to all conditional Poisson sampling schemes
satisfying the integrality assumption (A2).
To see how Theorem \ref{thm:ex.mle} differs from Theorem \ref{thm:mle},
a direct application of the latter yields that the MLE exists if and
only if $\bft$ belongs to the interior of the $(d-m)$-dimensional
polyhedron
\[
C_{\m{V}} : = \operatorname{conv}  ( \{ \bft\dvtx \bft= \m{A}^\top
\bfx, \bfx\in\mathbb{N}^\mathcal{I}, \m{V}^\top\bfx= \mathbf{1} \}  ).
\]
For Poisson sampling, this polyhedron is in fact the marginal cone,
and, for multinomial sampling,\vadjust{\goodbreak} it is the polytope $\{ \m{V} \mathbf{x}
\dvtx \mathbf{x} \in\operatorname{conv}(\m{A}) \}.$ Under product
multinomial sampling, $C_{\m{V}}$ is the Minkowsoki addition [see,
e.g., \citet{ZIE98}, \citet{SCH98}] of $m$ polytopes, one for each
multinomial, while under Poisson--multinomial scheme it is the Minkowski
sum of a polyhedral cone and as many polytopes as multinomial
constraints. Even though it has smaller ambient dimension than the
marginal cone, $C_{\m{V}}$ is a geometric object that can be rather
difficult to handle, both computationally and theoretically. In
contrast, we show that, for any sampling scheme satisfying conditions
(A0)--(A2), it is in fact sufficient to deal with the polyhedral cone
$C_{\m{A}}$, which is simpler to describe and analyze, both
algorithmically and in theory; see the supplementary material
\citet{FRsupp}. In \citet{NET} we provide various examples of how
Theorem~\ref{thm:ex.mle} can be used to simplify the task of
characterizing existence of the MLE for otherwise complicated models
for networks and random graphs. These particular models are based on
product multinomial sampling constraints, in which case Theorem
\ref{thm:ex.mle} yields what is known in polyhedral geometry as the
Cayley trick.


\subsection{Parameter estimability}\label{sec:estimab}

We now turn to the issue of estimability of the natural and log-linear
parameters when the MLE does not exist.
In our analysis, we rely on the key notion of facial sets, originally
introduced in a~slightly different form by \citet{GMS06}.

\begin{definition}\label{def:facial}
For a log-linear subspace $\mathcal{M}$, a set $\mathcal{F} \subseteq
\mathcal{I}$ is a facial set of $\mathcal{M}$, when, for some
$\bolds{\mu} \in\mathcal{M}$,
\begin{eqnarray*}
\bolds{\mu}(i) &=& 0 \qquad\mbox{if } i \in\mathcal{F},\\
\bolds{\mu}(i) &<& 0 \qquad\mbox{if } i \notin\mathcal{F}.
\end{eqnarray*}
\end{definition}

Equivalently, $\mathcal{F}$ is a facial set of $\mathcal{M}$ when, for
any design matrix $\m{A}$ for $\mathcal{M}$ (not necessarily of full
column rank), there exists some $\mathbf{c} \in\mathbb{R}^d$ such that
\begin{eqnarray}\label{eq:facial}
(\mathbf{a}_i,\mathbf{c}) &= &0 \qquad \mbox{if }  i \in\mathcal{F},\nonumber\\[-8pt]\\[-8pt]
(\mathbf{a}_i,\mathbf{c}) &< &0 \qquad\mbox{if }  i
\notin\mathcal{F},\nonumber
\end{eqnarray}
where $\mathbf{a}_i$ denotes the $i$th row of $\m{A}$. Facial sets
encode combinatorial and geometric properties of the log-linear
subspace $\mathcal{M}$ which turn out to be crucial to our analysis. We
summarize these properties in the next lemma.

\begin{lemma}\label{lem:F}
Let $\m{A}$ be a design matrix of $\mathcal{M}$. The lattice of facial
sets of~$\mathcal{M}$ is isomorphic to the face lattice of the marginal
cone $C_{\m{A}}$. In particular, $\mathcal{F}$ is a~facial set of
$\mathcal{M}$ if and only if $\{ \mathbf{a}_i, i \in \mathcal{F}\}$
span the face of $C_\m{A}$ isomorphic to $\mathcal{F}$.
\end{lemma}

Using this result, we can paraphrase Theorem \ref{thm:ex.mle} as
follows [compare with Theorem 3.2 in \citet{HAB74}]:\vadjust{\goodbreak}

\begin{corollary}
The MLE exists if and only if $\operatorname{supp}(\mathbf{n})$ is not
contained in any facial set of $\mathcal{M}$.
\end{corollary}

We describe algorithms for determining facial sets and for using the
previous corollary in \citet{FRsupp}.

\subsubsection{Estimability of the natural parameters}\label{sec:est.cond.poiss}

In this section, we rely on arguments proposed in \citet{ERGM09} to
study the estimability of the natural parameters. Let $C_{\m{V}}$
denote the convex support of the family arising from a conditional
Poisson scheme specified by a constraint matrix $\m{V}$; see Appendix
\ref{app:ext.expo}. Suppose that the observed sufficient statistics
$\mathbf{t} = \m{A}^\top\mathbf{n}$ belong to the relative interior of
face $F_{\m{V}}$ of $C_{\m{V}}$ of dimension $d_F$. Thus, the MLE of
the natural parameters for the original family, supported on
$S(\m{V})$, is nonexistent, but the MLE of the natural parameter of the
extended family supported $F_{\m{V}}$ is well defined.
Theorem \ref{lem:ZL} below generalizes Lemma \ref{lem:thetaN} by
showing that, when the MLE does not exist, the linear combinations of
the natural parameters that are estimable are determined, not only by
the deterministic linear subspace arising from the sampling
constraints, but also by the random linear subspace spanned by the
normal cone to the face $F$ of the marginal cone $C_{\m{A}}$ containing
$\m{A}^\top\mathbf{n}$ in its relative interior. As for the log-linear
parameter, nonexistence of the MLE entails that only points in
$\pi_{\mathcal{F}}  ( \mathcal{M} \ominus\mathcal{N}   )$ are
estimable, where $\mathcal{F}$ is the random facial set corresponding
to $F$.

In preparation for the result, we need to set up some additional
notation. By Lemma \ref{ref:homo} in Appendix \ref{appendix.proof},
there exists one face $F$ of $C_{\m{A}}$ of dimension $m + d_F$ that
contains $F_{\m{V}}$, with facial set $\mathcal{F}$. Let $N_F$ be the
normal cone to~$F$ and $\mathcal{L}_F \subset\mathbb{R}^d$ be the
linear subspace spanned by $N_F$, so that
$\operatorname{dim}(\mathcal{L}_F) = d - m - d_F$ (recall that, without
loss of generality, we assume $C_{\m{A}}$ to be full-dimensional). We
further define the linear subspace
\[
\mathcal{N}_F := \{ \m{A} \bolds{\beta}, \bolds{\beta} \in\mathcal {Z}
+ \mathcal{L}_F\},
\]
where $\mathcal{Z}$ is given in \eqref{eq:Z}. Just like in Lemma
\ref{lem:thetaN}, we define the following equivalence relation on
$\mathbb{R}^{d}$: $\theta_1 \stackrel{\mathcal{N}_F}{\sim} \theta
_2$ if
and only if $\bftheta_1 - \bftheta_2 \in\mathcal{Z} + L_F$, and write~$\bftheta_{\mathcal{N}_F}$ for the equivalence class containing
$\bftheta$. Finally, $\Theta_{\mathcal{N}_F} : = \{
\bftheta_{\mathcal{N}_F}, \bftheta\in\mathbb{R}^d\}$.\looseness=-1

\begin{theorem}\label{lem:ZL}
Consider the exponential family \eqref{eq:exp.A}, with $\m{A}$ of full
rank $d$, and suppose that conditions \textup{(A0)--(A2)} hold.  Let
$F_{\m{V}}$ be a face of the convex support and $\mathcal{F}$ the
corresponding facial set of the normal cone.
\begin{enumerate}[(ii)]
\item[(i)] For any
$\bftheta \in \mathbb{R}^d$, the set $\bftheta_{\mathcal{N}_F}$ is an
affine subspace of $\mathbb{R}^d$ of dimension $m +
\operatorname{dim}(\mathcal{L}_F) = d-d_F$. The set $\Theta_{\mathcal{N}_F}$
is a $d_F$-dimensional dimensional vector space isomorphic to
$\pi_{\mathcal{F}} ( \mathcal{M} \ominus \mathcal{N} )$ and is
comprised of parallel $(d - d_F)$-dimensional affine subspaces of
$\mathbb{R}^d$.
\item[(ii)] The extended family corresponding to
$F_{\m{V}}$  is non-identifiable: any two points $\bftheta_1
\stackrel{\mathcal{N}_F}{\sim}\bftheta_2$ specify the same
distribution. In fact, the family is parametrized by
$\Theta_{\mathcal{N}_F}$, or, equivalently, by $\pi_{\mathcal{F}} (
\mathcal{M} \ominus \mathcal{N} )$. Therefore, it is of order $d_F$.
\end{enumerate}
\end{theorem}

The main point of Theorem \ref{lem:ZL} is that only natural parameters
in $\Theta_{\mathcal{N}_F}$ [or the log-linear parameters in
$\pi_{\mathcal{F}}  ( \mathcal{M} \ominus\mathcal{N}   )$] are
estimable, with both sets being now random. In principle,
nonidentifiability of the natural parameters, due to a nonexistent MLE,
can be resolved using the same procedure of reduction to minimality
described in the remarks following Lemma \ref{lem:thetaN}: identify a
set of linearly independent vectors in $\mathbb{R}^{\mathcal{I}}$
spanning $\mathcal{M} \cap \mathcal{N}^\bot_F$, and use them to build a
new design matrix of dimension $I \times d_F$. However, unlike the
reduction to minimality carried out to remove the effect of the
sampling constraints, which is design-dependent but not data-dependent,
this reduction depends on the random subspace $\mathcal{N}_F$ (the
randomness arising from the exposed face $F$). Furthermore, while the
sampling constraint reduction is easy to implement, since the matrix
$\m{V}$ is known, this second reduction requires us to compute a basis
for $\mathcal{L}_F$, the linear space spanned by the normal cone to
$F$. For the mean value parameter, the problem is to compute the facial
set associated to the face $F$ based solely on the observed sufficient
statistics $\mathbf{t}$, which amounts to identifying the face of
$C_{\m{A}}$ containing $\mathbf{t}$ in its relative interior. In
general, both of these tasks are highly nontrivial, due to the
combinatorial complexity of the face lattice of $C_{\m{A}}$; see the
examples in Section \ref{sec:inf}. In the supplementary material
[\citet{FRsupp}], we describe algorithms for accomplishing these tasks.

As a corollary to Theorem \ref{lem:ZL}, we can obtain each family in
the extended family via a conditional Poisson sampling scheme that
forces the base measure to be supported on $F_{\m{V}}$, or
equivalently, by requiring that the cells in $\mathcal{F}^c$ have zero
probability of containing positive counts. In this case, it is clear
that assumption (A1) is violated. As a result, we can view each such
family as a log-linear model under Poisson sampling scheme containing
structural zeros along the (random) coordinates $\mathcal{F}^c$. This
is in fact consistent with the interpretation by \citet{BARN78}, page
156, of the extended MLE as a conditional MLE, given that sufficient
statistics lie on the boundary of the convex support. We formalize this
observation in the next result.

\begin{corollary}\label{cor:struct}
Each face $F$ of $C_{\m{V}}$ of dimension $0 \leq d_F \leq d - m$ can
be obtained as the convex support corresponding to the conditional
Poisson scheme with constraint subspace $\mathcal{N}_{F}$, where
$\operatorname{dim}(\mathcal{N}_{F}) = d - d_F$.
\end{corollary}

Using the same arguments as in the remarks following Lemma
\ref{lem:thetaN}, we also see that the Fisher information matrix at the
extended MLE has rank $d_F < d$, and therefore, is rank-deficient. This
remains the case, even after accounting for\vadjust{\goodbreak} the sampling constraints.
Statistically, the singularity of the observed Fisher information
implies that the standard errors are not defined. From an algorithmic
standpoint, this observation implies that the Newton--Raphson method
for computing the MLE is bound to run into numerical instabilities, due
to the fact that the Hessian matrix of the log-likelihood function is
singular at any optimum [an issue illustrated empirically in
\citet{ROY05}]. Furthermore, Corollary~2.8 in \citet{ERGM09} shows
that, under a nonexistent MLE, every point in the normal cone $N_F$ to
the face $F$ containing the observed sufficient statistics is a
(random) direction of recession of the negative log-likelihood
function, so that there are infinitely many directions of maximal
increase of the log-likelihood function.


\subsubsection{Estimability of the mean value parameters under Poisson
and product multinomial schemes}\label{sec:mean.value}

We now specialize our analysis to the case of Poisson and product
multinomial sampling schemes. Besides their popularity, the main reason
for focusing on these two particular sampling schemes is that the
estimates of the cell mean values are highly interpretable. Under the
Poisson scheme, the cell mean values are just the expected cell counts,
while under the product multinomial scheme they are the conditional
expectations of the cell counts given the grand total (in the
multinomial case) or given the total counts in the portions of the
table associated with the partitions used to define the product
multinomial constraints. For other conditional Poisson sampling
schemes, not only are the conditional cell mean values difficult to
compute due to the unknown normalizing constant, but they are also less
interpretable.


Following \citet{LAU96}, we consider $\overline{M} = \operatorname{cl} (
\{ \exp(\bolds{\mu}), \bolds{\mu} \in\mathcal{M} \}   )$, the closure
of the set of all cell mean values for a log-linear subspace
$\mathcal{M}$. Thus, $\mathbf{m} \in\overline{M}$ if and only if
$\mathbf{m} = \lim_n e^{\bolds{\mu}_n}$, for some sequence $\{
\bolds{\mu}_n \}_n \subset\mathcal{M}$. \citet{LAU96} calls the set
$\overline{M}$ the {\it extended log-affine model}.

\begin{theorem}\label{thm:main}
Let $\mathbf{t}$ be the observed sufficient statistics, and let
$\mathcal{F}$ be facial set corresponding to the face of $C_{\m{A}}$
containing $\mathbf{t}$ in its relative interior. The MLE of the cell
mean vector exists, is unique and identical under Poisson and product
multinomial if and only if $\mathcal{F} = \mathcal{I}$. If $\mathcal{F}
\subsetneq\mathcal{I}$, there exists one point
$\widehat{\mathbf{m}}{}^{\m{e}}$ in $\overline{M}$ such that
$\widehat{\mathbf{m}}{}^{\m{e}} = \lim_n \exp(\bolds{\mu}_n )$, where
$\{\bolds{\mu}_n \}_n \subset\mathcal{M}$ is any optimizing sequence
such that
\[
\lim_n \ell^P(\bolds{\mu}_n) = \sup_{\bolds{\mu} \in\mathcal{M}}
\ell^P(\mu)  \quad\mbox{and}\quad  \lim_n \ell^M(\bolds{\mu}_n) =
\sup_{\tilde{\bolds{\mu}} \in\widetilde{\mathcal{M}}} \ell
^M(\tilde{\bolds{\mu}}).
\]
Furthermore, $\operatorname{supp}(\widehat{\mathbf{m}}{}^{\m{e}}) =
\mathcal{F}$ and $\Pi_{\mathcal{M}} \mathbf{n} = \Pi_{\mathcal{M}}
\widehat{\mathbf{m}}{}^{\m{e}}$.
\end{theorem}

This result shows that, for any observed table $\mathbf{n}$, the
log-likelihood functions in both sampling schemes admits always a
unique maximizer, $\widehat{\mathbf{m}}{}^{\m{e}}$. Though
supported
only on the facial set associated with $\mathbf{t}$, this vector
exhibits exactly\vadjust{\goodbreak} the same features as the ``ordinary'' MLE: it is the
unique point $\widehat{\mathbf{m}}{}^{\m{e}} \in\overline{M} $ such
that $\m{A}^\top\widehat{\mathbf{m}}{}^{\m{e}} = \m{A}^\top\mathbf{n} $
and provided that $\mathcal{N} \subset\mathcal{M}$, maximizes both the
Poisson and product multinomial likelihoods. The substantial difference
is that $\widehat{\mathbf{m}}{}^{\m{e}}$ has positive coordinates only
along the cells in the facial set $\mathcal{F}$. Theorem \ref{thm:main}
generalizes Theorem 4.8 in \citet{LAU96}. The improvement consists of
identifying exactly the supports of the limit points in $\overline{M}$,
which are precisely the facial sets of $C_{\m{A}}$.\vspace*{-2pt}

\begin{definition}
The vector $\widehat{\mathbf{m}}{}^{\m{e}}$ is the \textit{extended
MLE} of $\mathbf{m}$ and the zeros appearing in along the coordinates
in $\mathcal{F}^c = \mathcal{I} \setminus\mathcal{F}$ are called the
{\it likelihood zeros}.\vspace*{-2pt}
\end{definition}

The term likelihood zeros highlight the fact that those zero counts,
though arising as sampling and not as structural zeros, have a
significant impact on the likelihood function and its optimizers.\vspace*{-2pt}

\subsection{The geometry of the extended Poisson family}\label{sec:variety}

The results of Theorem \ref{thm:main} suggest that, for the Poisson and
product multinomial schemes, we could, in fact, take the set
$\overline{M}$ to be the cell mean value parameter space for the
extended exponential family of distributions for the actual contingency
table, not its sufficient statistics. We formalize this idea by relying
on geometric considerations. For ease of readability, and without loss
of generality, we focus on the Poisson sampling scheme, and only
sketch how our results apply also to product multinomial cases.

For a vector $\mathbf{u} \in\mathbb{R}^\mathcal{I}$, let
\[
\mathbf{u}^+ = \bigl\{ \max\{ \mathbf{u}(i),0 \}, i \in\mathcal{I}\bigr\}
\quad\mbox{and}\quad  \mathbf{u}^- = \bigl\{ \min\{ \mathbf{u}(i),0 \}, i
\in \mathcal{I}\bigr\},
\]
so that $\mathbf{u} = \mathbf{u}^+ - \mathbf{u}^-$ and
$\operatorname{supp}(\mathbf{u}^+)
\cap\operatorname{supp}(\mathbf{u}^-) =\varnothing$. Furthermore, for any
pair of nonnegative vectors $\mathbf{x}$ and $\mathbf{u}$ in
$\mathbb{R}^\mathcal{I}$, write
\[
\mathbf{x}^\mathbf{u} = \prod_{i} \mathbf{x}(i)^{\mathbf{u}(i)}
\]
for the associated monomial. Following \citet{GMS06}, page~1469 and
Lemma A.1, we consider the toric variety $X_\mathcal{M}$ corresponding
to the log-linear model $\mathcal{M}$.\vspace*{-2pt}

\begin{definition}
The nonnegative toric variety $X_\mathcal{M}$ associated to the
log-linear subspace $\mathcal{M}$ is the set of all vectors $\mathbf{x}
\in\mathbb{R}^\mathcal{I}_{\geq0}$ such that
\begin{equation}\label{eq:toric}
\mathbf{x}^{\mathbf{u}^+} = \mathbf{x}^{\mathbf{u}^-} \qquad
 \forall\mathbf{u} \in \mathcal{M}^\bot.\vspace*{-2pt}
\end{equation}
\end{definition}

Geometrically, $X_\mathcal{M}$ is the intersection of the solution set
of a system of polynomial equations with the nonnegative orthant. It is
easy to see that any $\mathbf{m}>0$ such that $\log(\mathbf{m}) \in
\mathcal{M}$ satisfies \eqref{eq:toric}. Equation \eqref{eq:toric} can
still hold, however, when some of the coordinates of $\mathbf{m}$ are
zero. Finally, for any $\bolds{\xi} \in C_{\m{A}}$, consider the
polyhedron
\begin{equation}
P_{\bolds{\xi}} = \{\mathbf{x} \in\mathbb{R}^{\mathcal{I}}_{\geq0}
\dvtx  \mathrm{A} \mathbf{x} = \bolds{\xi} \}.\vadjust{\goodbreak}
\end{equation}
For a given sufficient statistic $\mathbf{t = \m{A} \mathbf{n}}$, the
set of lattice points in $P_{\mathbf{t}}$, known as the {\it fiber} of
$\mathbf{t}$, consists of all possible tables having the same
sufficient statistics as the observed table $\mathbf{n}$.\vspace*{-2pt}

\begin{theorem}\label{thm:geom}
\textup{(i)} $\overline{M} = X_{\mathcal{M}}$.\vspace*{-6pt}
\begin{enumerate}[(iii)]
\item[(ii)] For any nonzero $\mathbf{m} \in X_{\mathcal {M}}$,
$\operatorname{supp}(\mathbf{m})$ is a facial set of $C_{\m{A}}$.

\item[(iii)] The linear map $\m{A} \dvtx
\mathbb{R}^{\mathcal{I}} \rightarrow\mathbb{R}^d$, given by $\mathbf{m}
\mapsto\m{A} \mathbf{m}$, defines a homeomorphism between
$X_{\mathcal{M}}$ and $C_{\m{A}}$.
\item[(iv)] For any observable sufficient statistic
$\mathbf{t} = \m{A}\mathbf{n},$ $\{ \widehat{\mathbf{m}}{}^{\m{e}} \} =
X_{\mathcal{M}} \cap P_{\mathbf{t}}$ and
$\widehat{\mathbf{m}}{}^{\m{e}} \in\operatorname{ri}(P_{\mathbf{t}})$.\vspace*{-2pt}
\end{enumerate}
\end{theorem}

Part (i) of Theorem \ref{thm:geom} is due to \citet{GMS06}, while a
slightly less general version of part (iii) is a standard result in the
algebraic statistics literature; see, for example, \citet{PS05},
\citet{DSS08}.


Overall, Theorem \ref{thm:geom} shows that the set $\overline{M} $ is
homeomorphic to the margin\-al cone $C_{\m{A}}$ and, therefore, as
anticipated, we can use it as a legitimate mean value space for the
extended family of the cell counts. The advantage of $\overline{M} $
over $C_{\m{A}}$ is its direct interpretability in terms of cell mean
values.
This result extends directly to the multinomial sampling scheme. In
this case, $\m{A}$ specifies a homeomorphism between $\{\mathbf{x} \in
X_{\mathcal{M}} \dvtx \sum_i \mathbf{x}_i = 1\} $ and $P_{\m{A}} =
\operatorname{conv}(\m{A})$, which is known in algebraic geometry as
the moment map; see \citet{FUL78}, \citet{EWA96}. In fact, under
multinomial scheme, the extended mean-value space can be taken to be
the intersection of $X_{\mathcal{M}}$ with the probability simplex in~$\mathbb{R}^{\mathcal{I}}$. Furthermore, since $\mathcal{M}$ contains
the constant functions, $P_{\m{A}}$ and~$C_{\m{A}}$ have identical
facial sets.
For product multinomial sampling schemes, a characterization of the
mean value space analogous to the one given in Theorem~\ref{thm:geom}
is also possible, though somewhat more involved. We refer the reader to
\citet{JASON08} for details and a different derivation. In this
particular case, the convex support arises as a Minkwoski sum of
polytopes, one for every multinomial. Then, the proof of Theorem~\ref{thm:ex.mle} reveals that facial sets of the convex support are
also facial sets of the marginal cone, even though the opposite is not
true. See \citet{NET} for an application of these results to network
models.

Finally, part  (iv) of Theorem \ref{thm:geom} shows that the extended
MLE is the only point in $P_{\mathbf{t}}$ satisfying the log-linear
model conditions. This result can be also interpreted in terms of
I-divergence projections [\citeauthor{CS75} (\citeyear{CS75,CSI89})], and provides
the geometric basis for showing convergence of iterative methods for
extended maximum likelihood estimation such as the iterative
proportional scaling algorithm of \citet{DARRAT72}. In the interest of
space, we do not pursue this analysis.\vspace*{-2pt}

\section{Inference under a nonexistent MLE}\label{sec:inf}

We have shown that when the MLE does not exist, only some of the model
parameters (both under the natural and mean-value parametrization) are
estimable, and we\vadjust{\goodbreak} have identified the parameters that can instead be
estimated within the extended family. Thus, when the MLE is
nonexistent, statistical inference is still feasible, but only for the
reduced family whose parameters are fully estimable.

As described at the end of Section \ref{sec:est.cond.poiss}, we can
obtain the relevant extended exponential family by computing a new
random design matrix $\m{A}_\mathcal{F}$ whose column span is
$\pi_{\mathcal{F}}(\mathcal{M} \ominus\mathcal{N})$, where
$\mathcal{F}$ is the random facial set corresponding to the face $F$ of
the marginal cone containing the sufficient statistics in its relative
interior. We can then use this new design matrix to specify a new
exponential family as in~\eqref{eq:exp.A}, where {\it only the cells
$\mathcal{F}$ have positive probability of being observed.} We carry
out inference within this extended family or, equivalently,
conditionally on the sufficient statistics being on the face~$F$, as
advocated by \citet{BARN78}, page~156. By Corollary \ref{cor:struct},
this is equivalent to treating the coordinates in $\mathcal{F}$ {\it as
if they were structural zeros.} Thus, dealing with a nonexistent MLE
reduces, in practice, to fitting the same log-linear model under the
additional (random) constraints that the cells in $\mathcal{F}^c$,
which are not estimable, be treated as structural zeros. The same
approach is also advocated in \citet{GEYER09}. In practice, this
entails replacing the MLE with the extended MLE and, quite importantly,
adjusting the number of degrees of freedom, now to be computed as the
difference between the cardinality of the facial set $|\mathcal{F}|$
(i.e., the number of cell mean values that can be estimated), and the
number of estimable parameters, namely $\operatorname{dim}  (
\pi_{\mathcal{F}}(\mathcal{M} \ominus\mathcal{N}) ) =
\operatorname{dim}(F) - m$. Using the adjusted number of degrees of
freedom, asymptotic $\chi^2$ tests for goodness of fit [see, e.g.,
\citet{CR88}] can then still be applied.
Algorithms for carrying out the numerical tasks just described are
presented in the supplementary material [\citet{FRsupp}].


\section{Examples of likelihood zeros}

Below, we illustrate by means of examples various practical aspects of
goodness-of-fit testing when the MLE is nonexistent, and we show how to
appropriately adjust the number of degrees of freedom. We will focus on
hierarchical log-linear models [see, e.g., \citet{BFH75}], and refer
the reader to \citet{DFRSZ09} and \citet{NET} for other examples of
this kind.

Our polyhedral characterization of the conditions for the existence of
the MLE permits to generate novel examples of patterns of sampling
zeros causing nonexistence of the MLE for hierarchical log-linear
models without producing null margins, an instance that is virtually
ignored in all statistical software. As pointed out by \citet{ROY05},
the {\tt R} [\citet{R}] routines {\tt loglin} and {\tt glm}, as well as
virtually any other software for inference and model selection for
log-linear models, does no detect nonexistence and report the
unadjusted, incorrect, numbers of degrees of freedom for all the
examples below. In the analysis of sparse tables, it is also common
practice to add small positive quantities to the zero cells, in order
to avoid numerical issues with the computation of the MLE. We remain
highly skeptical of the numerical advantages of this ad-hoc procedure,
and remark that such adjustments will make it impossible to detect
nonexistence of the MLE and to distinguish the estimable parameters.



The examples of likelihood zeros in Examples \ref{ex2}--\ref{ex4}
suggest that the combinatorial complexity of hierarchical log-linear
models, measured by the number of facets of the marginal cone, can be
quite significant. In the reported examples, as well as in many other
experiments we conducted, for many models the number of facets
associated with zero margins appears to be much smaller than the total
number of facets, indicating that, at least combinatorially, likelihood
zeros associated to positive margins are much more frequent (though
never detected). Below we use the classic notation to represent the
generating class of a hierarchical log-linear model; for example, see
\citet{BFH75}. Empty cells indicate positive counts. All the
calculations were carried out in {\tt polymake}
[\citet{polymake}].\looseness=-1


\begin{example}\label{ex1} The $2^3$ table and the model $[12][13][23]$ of
no-second-order interaction. The MLE is not defined because the two
likelihood zeros expose one of the 16 facets of the marginal cone. This
example, due to \citet{HAB74}, was the only published example a
log-linear model with nonexistent MLE and positive margins; see
\citet{ROY05}, Section 5, for a~general result concerning binary
$K$-way tables and the model of no-$(K-1)$st interaction\vspace*{7pt}.
\begin{center}
%
\begin{tabular}{cc}
\hline\\[-13pt]
\multicolumn{1}{|c|}{0} &\multicolumn{1}{c@{}|}{\hphantom{0\hspace*{5.8pt}}} \\[-1pt]
\hline\\[-13pt]
\multicolumn{1}{|c|}{} &\multicolumn{1}{c@{}|}{} \\
\hline
\end{tabular}\hspace*{10pt}
\begin{tabular}{cc}
\hline\\[-13pt]
\multicolumn{1}{|c|}{\hphantom{0}} &\multicolumn{1}{c@{}|}{\hphantom{0\hspace*{5.8pt}}} \\[-1pt]
\hline\\[-13pt]
\multicolumn{1}{|c|}{\hphantom{0}}&\multicolumn{1}{@{}c@{}|}{0}
\\
\hline
\end{tabular}
%
\end{center}\vspace*{7pt}
The dimension of the log-linear subspace for this model, or,
equivalently, of the marginal cone, is $7$, leaving $1$ degree of
freedom when the MLE exists. However, because of the likelihood zeros,
inference can only be made for the $6$-dimensional exposed facet. Since
the cardinality of the associated facial set~$\mathcal{F}$ is also $6$,
the resulting extended log-linear model is the saturated model on
$\mathcal{F}$.
\end{example}


\begin{example}\label{ex2} The $3^3$ table and the model
$[12][13][23]$. The MLE is not defined because the pattern of
likelihood zeros exposes one of the 207 facets of the marginal cone. Of
all the facets, only 27 are associated to zero margins.\vspace*{7pt}
%
\begin{center}
%
%
\begin{tabular}{ccc}
\hline\\[-13pt]
\multicolumn{1}{|c|}{0} &\multicolumn{1}{c}{}&\multicolumn{1}{|c@{}|}{\hphantom{0\hspace*{5.8pt}}} \\[-1pt]
\hline\\[-13pt]
\multicolumn{1}{|c|}{} &\multicolumn{1}{c}{}&\multicolumn{1}{|c@{}|}{\hphantom{0\hspace*{5.8pt}}} \\[-1pt]
\hline\\[-13pt]
\multicolumn{1}{|c|}{0} &\multicolumn{1}{c}{0}&\multicolumn{1}{|c@{}|}{\hphantom{0\hspace*{5.8pt}}} \\
\hline
\end{tabular}\hspace*{10pt}
\begin{tabular}{ccc}
\hline\\[-13pt]
\multicolumn{1}{|c|}{\hphantom{0}} &\multicolumn{1}{c}{}&\multicolumn{1}{|c@{}|}{\hphantom{0\hspace*{5.8pt}}} \\[-1pt]
\hline\\[-13pt]
\multicolumn{1}{|c|}{} &\multicolumn{1}{c}{0}&\multicolumn{1}{|c@{}|}{0} \\[-1pt]
\hline\\[-13pt]
\multicolumn{1}{|c|}{} &\multicolumn{1}{c}{0}&\multicolumn{1}{|c@{}|}{} \\
\hline
\end{tabular}\hspace*{10pt}
\begin{tabular}{ccc}
\hline\\[-13pt]
\multicolumn{1}{|c|}{0} &\multicolumn{1}{c}{\hphantom{0}}&\multicolumn{1}{|@{}c@{}|}{0} \\[-1pt]
\hline\\[-13pt]
\multicolumn{1}{|c|}{} &\multicolumn{1}{c}{}&\multicolumn{1}{|@{}c@{}|}{0} \\[-1pt]
\hline\\[-13pt]
\multicolumn{1}{|c|}{} &\multicolumn{1}{c}{}&\multicolumn{1}{|c@{}|}{\hphantom{0\hspace*{5.8pt}}} \\
\hline
\end{tabular}  
\end{center}\vspace*{7pt}
The dimension of the facet is $18$, which is also the cardinality of
the facial set for this configuration of likelihood zeros. As in
the
previous example, this defines the saturated model on $\mathcal{F}$,
giving $0$ adjusted degrees of freedom and making $\chi^2$
approximations not applicable.\vadjust{\goodbreak}

Under the same log-linear model, the MLE does not exist also when the
following pattern of zeros arises:\vspace*{7pt}
\begin{center}
\begin{tabular}{ccc}
\hline\\[-13pt]
\multicolumn{1}{|c|}{0} &\multicolumn{1}{c}{}&\multicolumn{1}{|c@{}|}{\hphantom{0\hspace*{5.8pt}}} \\[-1pt]
\hline\\[-13pt]
\multicolumn{1}{|c|}{0} &\multicolumn{1}{c}{}&\multicolumn{1}{|c@{}|}{\hphantom{0\hspace*{5.8pt}}} \\[-1pt]
\hline\\[-13pt]
\multicolumn{1}{|c|}{} &\multicolumn{1}{c}{\hphantom{0}}&\multicolumn{1}{|c@{}|}{\hphantom{0\hspace*{5.8pt}}} \\
\hline
\end{tabular}\hspace*{10pt}
\begin{tabular}{ccc}
\hline\\[-13pt]
\multicolumn{1}{|c|}{\hphantom{0}} &\multicolumn{1}{c}{}&\multicolumn{1}{|c@{}|}{\hphantom{0\hspace*{5.8pt}}} \\[-1pt]
\hline\\[-13pt]
\multicolumn{1}{|c|}{} &\multicolumn{1}{c}{}&\multicolumn{1}{|@{}c@{}|}{0} \\[-1pt]
\hline\\[-13pt]
\multicolumn{1}{|c|}{} &\multicolumn{1}{c}{\textbf{0}}&\multicolumn{1}{|@{}c@{}|}{0} \\
\hline
\end{tabular}\hspace*{10pt}
\begin{tabular}{ccc}
\hline\\[-13pt]
\multicolumn{1}{|c|}{0} &\multicolumn{1}{c}{\hphantom{0}}&\multicolumn{1}{|@{}c@{}|}{\textbf{0}} \\[-1pt]
\hline\\[-13pt]
\multicolumn{1}{|c|}{0} &\multicolumn{1}{c}{}&\multicolumn{1}{|c@{}|}{\hphantom{0\hspace*{5.8pt}}} \\[-1pt]
\hline\\[-13pt]
\multicolumn{1}{|c|}{} &\multicolumn{1}{c}{}&\multicolumn{1}{|c@{}|}{} \\
\hline
\end{tabular}
\end{center}\vspace*{7pt}
In this example, the zeros displayed in bold are not likelihood zeros,
but the others are. Indeed, their presence or absence has no effect on
the existence of the MLE. Furthermore, when the extended MLE is
computed, the boldfaced zero counts will be replaced by positive
entries, while the likelihood zeros will stay zero. The number of
degrees of freedom in this example is $3$, because the total number of
estimable cell mean values is $21$, and the number of parameters for
the reduced model is $18$.

In our last example, the MLE is defined, despite the table being very
sparse, because no facet of the marginal cone is exposed [source:
\citet{ROY05}].\vspace*{7pt}
\begin{center}
\begin{tabular}{ccc}
\hline\\[-13pt]
\multicolumn{1}{|c|}{} &\multicolumn{1}{c}{0}&\multicolumn{1}{|@{}c@{}|}{0} \\[-1pt]
\hline\\[-13pt]
\multicolumn{1}{|c|}{0} &\multicolumn{1}{c}{}&\multicolumn{1}{|@{}c@{}|}{0} \\[-1pt]
\hline\\[-13pt]
\multicolumn{1}{|c|}{0} &\multicolumn{1}{c}{0}&\multicolumn{1}{|c@{}|}{\hphantom{0\hspace*{5.8pt}}} \\
\hline
\end{tabular}\hspace*{10pt}
\begin{tabular}{ccc}
\hline\\[-13pt]
\multicolumn{1}{|c|}{0} &\multicolumn{1}{c}{0}&\multicolumn{1}{|c@{}|}{\hphantom{0\hspace*{5.8pt}}} \\[-1pt]
\hline\\[-13pt]
\multicolumn{1}{|c|}{} &\multicolumn{1}{c}{0}&\multicolumn{1}{|@{}c@{}|}{0} \\[-1pt]
\hline\\[-13pt]
\multicolumn{1}{|c|}{0} &\multicolumn{1}{c}{}&\multicolumn{1}{|@{}c@{}|}{0} \\
\hline
\end{tabular}\hspace*{10pt}
\begin{tabular}{ccc}
\hline\\[-13pt]
\multicolumn{1}{|c|}{0} &\multicolumn{1}{c}{}&\multicolumn{1}{|@{}c@{}|}{0} \\[-1pt]
\hline\\[-13pt]
\multicolumn{1}{|c|}{0} &\multicolumn{1}{c}{0}&\multicolumn{1}{|c@{}|}{\hphantom{0\hspace*{5.8pt}}} \\[-1pt]
\hline\\[-13pt]
\multicolumn{1}{|c|}{} &\multicolumn{1}{c}{0}&\multicolumn{1}{|@{}c@{}|}{0} \\
\hline
\end{tabular}
\end{center}
\end{example}

\begin{example}\label{ex3} The $4 \times4 \times4$ table and the
model $[12][13][23]$. The MLE is not defined because the pattern of
zeros exposes one of 113,740 facets of the marginal cone [source:
\citet{MLE06}]. Of these, only 48 are associated to zero
margins.\vspace*{7pt}
\begin{center}
\begin{tabular}{cccc}
\hline\\[-13pt]
\multicolumn{1}{|c|}{0}  &   \multicolumn{1}{c|}{0} &   \multicolumn{1}{c}{0}  &    \multicolumn{1}{|@{}c@{}|}{} \\[-1pt]
\hline\\[-13pt]
\multicolumn{1}{|c|}{0} &   \multicolumn{1}{c|}{0} &   \multicolumn{1}{c}{}   &    \multicolumn{1}{|@{}c@{}|}{} \\[-1pt]
\hline\\[-13pt]
\multicolumn{1}{|c|}{0} &   \multicolumn{1}{c|}{} &   \multicolumn{1}{c}{}  &    \multicolumn{1}{|c@{}|}{\hphantom{0\hspace*{5.8pt}}} \\[-1pt]
\hline\\[-13pt]
\multicolumn{1}{|c|}{}  &   \multicolumn{1}{c|}{} &   \multicolumn{1}{c}{}  &    \multicolumn{1}{|c@{}|}{\hphantom{0\hspace*{5.8pt}}} \\
\hline
\end{tabular}\hspace*{10pt}
\begin{tabular}{cccc}
\hline\\[-13pt]
\multicolumn{1}{|c|}{\hphantom{0}}  &   \multicolumn{1}{c|}{0} &   \multicolumn{1}{c}{0}  &    \multicolumn{1}{|@{}c@{}|}{} \\[-1pt]
\hline\\[-13pt]
\multicolumn{1}{|c|}{} &   \multicolumn{1}{c|}{0} &   \multicolumn{1}{c}{}   &    \multicolumn{1}{|@{}c@{}|}{} \\[-1pt]
\hline\\[-13pt]
\multicolumn{1}{|c|}{} &   \multicolumn{1}{c|}{} &   \multicolumn{1}{c}{}  &    \multicolumn{1}{|c@{}|}{\hphantom{0\hspace*{5.8pt}}} \\[-1pt]
\hline\\[-13pt]
\multicolumn{1}{|c|}{} &   \multicolumn{1}{c|}{0} &   \multicolumn{1}{c}{0}  &    \multicolumn{1}{|@{}c@{}|}{0} \\
\hline
\end{tabular}\hspace*{10pt}
\begin{tabular}{cccc}
\hline\\[-13pt]
\multicolumn{1}{|c|}{}  &   \multicolumn{1}{c|}{\hphantom{0}} &   \multicolumn{1}{c}{0}  &    \multicolumn{1}{|@{}c@{}|}{} \\[-1pt]
\hline\\[-13pt]
\multicolumn{1}{|c|}{} &   \multicolumn{1}{c|}{} &   \multicolumn{1}{c}{}   &    \multicolumn{1}{|c@{}|}{\hphantom{0\hspace*{5.8pt}}} \\[-1pt]
\hline\\[-13pt]
\multicolumn{1}{|c|}{0} &   \multicolumn{1}{c|}{} &   \multicolumn{1}{c}{0}  &    \multicolumn{1}{|@{}c@{}|}{0} \\[-1pt]
\hline\\[-13pt]
\multicolumn{1}{|c|}{} &   \multicolumn{1}{c|}{} &   \multicolumn{1}{c}{0}  &    \multicolumn{1}{|@{}c@{}|}{0} \\
\hline
\end{tabular}\hspace*{10pt}
\begin{tabular}{cccc}
\hline\\[-13pt]
\multicolumn{1}{|c|}{}  &   \multicolumn{1}{c|}{} &   \multicolumn{1}{c}{\hphantom{0}}  &    \multicolumn{1}{|c@{}|}{\hphantom{0\hspace*{5.8pt}}} \\[-1pt]
\hline\\[-13pt]
\multicolumn{1}{|c|}{0} &   \multicolumn{1}{c|}{0} &   \multicolumn{1}{c}{}   &    \multicolumn{1}{|@{}c@{}|}{0} \\[-1pt]
\hline\\[-13pt]
\multicolumn{1}{|c|}{0} &   \multicolumn{1}{c|}{} &   \multicolumn{1}{c}{}  &    \multicolumn{1}{|@{}c@{}|}{0} \\[-1pt]
\hline\\[-13pt]
\multicolumn{1}{|c|}{} &   \multicolumn{1}{c|}{} &   \multicolumn{1}{c}{}  &    \multicolumn{1}{|@{}c@{}|}{0} \\
\hline
\end{tabular}
\end{center}
\end{example}

\begin{example}\label{ex4} The $3^4$ table and the 4-cycle model
$[12][14][23][34]$. The MLE is not defined because the pattern of zeros
exposes one of the 1116 facets of the marginal cone. Of these, only 36
are associated to zero margins.\vspace*{7pt}
\begin{center}
\begin{tabular}{ccc}
\hline\\[-13pt]
\multicolumn{1}{|c|}{0} &\multicolumn{1}{c}{\hphantom{0}}&\multicolumn{1}{|@{}c@{}|}{0} \\[-1pt]
\hline\\[-13pt]
\multicolumn{1}{|c|}{} &\multicolumn{1}{c}{}&\multicolumn{1}{|@{}c@{}|}{0} \\[-1pt]
\hline\\[-13pt]
\multicolumn{1}{|c|}{} &\multicolumn{1}{c}{}&\multicolumn{1}{|c@{}|}{\hphantom{0\hspace*{5.8pt}}} \\
\hline
\end{tabular}\hspace*{10pt}
\begin{tabular}{ccc}
\hline\\[-13pt]
\multicolumn{1}{|c|}{0} &\multicolumn{1}{c}{0}&\multicolumn{1}{|@{}c@{}|}{0} \\[-1pt]
\hline\\[-13pt]
\multicolumn{1}{|c|}{0} &\multicolumn{1}{c}{0}&\multicolumn{1}{|@{}c@{}|}{0} \\[-1pt]
\hline\\[-13pt]
\multicolumn{1}{|c|}{0} &\multicolumn{1}{c}{0}&\multicolumn{1}{|c@{}|}{\hphantom{0\hspace*{5.8pt}}} \\
\hline
\end{tabular}\hspace*{10pt}
\begin{tabular}{ccc}
\hline\\[-13pt]
\multicolumn{1}{|c|}{0} &\multicolumn{1}{c}{0}&\multicolumn{1}{|@{}c@{}|}{0} \\[-1pt]
\hline\\[-13pt]
\multicolumn{1}{|c|}{} &\multicolumn{1}{c}{0}&\multicolumn{1}{|@{}c@{}|}{0} \\[-1pt]
\hline\\[-13pt]
\multicolumn{1}{|c|}{} &\multicolumn{1}{c}{0}&\multicolumn{1}{|c@{}|}{\hphantom{0\hspace*{5.8pt}}} \\
\hline
\end{tabular}\vspace*{6pt}
\end{center}
\begin{center}
\begin{tabular}{ccc}
\hline\\[-13pt]
\multicolumn{1}{|c|}{0} &\multicolumn{1}{c}{}&\multicolumn{1}{|@{}c@{}|}{\hspace*{6pt}0\hspace*{6pt}} \\[-1pt]
\hline\\[-13pt]
\multicolumn{1}{|c|}{0} &\multicolumn{1}{c}{0}&\multicolumn{1}{|@{}c@{}|}{0} \\[-1pt]
\hline\\[-13pt]
\multicolumn{1}{|c|}{0} &\multicolumn{1}{c}{0}&\multicolumn{1}{|@{}c@{}|}{0} \\
\hline
\end{tabular}\hspace*{10pt}
\begin{tabular}{ccc}
\hline\\[-13pt]
\multicolumn{1}{|c|}{0} &\multicolumn{1}{c}{}&\multicolumn{1}{|c@{}|}{\hphantom{0\hspace*{5.8pt}}} \\[-1pt]
\hline\\[-13pt]
\multicolumn{1}{|c|}{0} &\multicolumn{1}{c}{0}&\multicolumn{1}{|@{}c@{}|}{0} \\[-1pt]
\hline\\[-13pt]
\multicolumn{1}{|c|}{0} &\multicolumn{1}{c}{0}&\multicolumn{1}{|@{}c@{}|}{} \\
\hline
\end{tabular}\hspace*{10pt}
\begin{tabular}{ccc}
\hline\\[-13pt]
\multicolumn{1}{|c|}{\hphantom{0}} &\multicolumn{1}{c}{}&\multicolumn{1}{|c@{}|}{\hphantom{0\hspace*{5.8pt}}} \\[-1pt]
\hline\\[-13pt]
\multicolumn{1}{|c|}{} &\multicolumn{1}{c}{0}&\multicolumn{1}{|@{}c@{}|}{0} \\[-1pt]
\hline\\[-13pt]
\multicolumn{1}{|c|}{} &\multicolumn{1}{c}{0}&\multicolumn{1}{|@{}c@{}|}{} \\
\hline
\end{tabular}\vspace*{6pt}
\end{center}
\begin{center}
\begin{tabular}{ccc}
\hline\\[-13pt]
\multicolumn{1}{|c|}{0} &\multicolumn{1}{c}{}&\multicolumn{1}{|@{}c@{}|}{\hspace*{6pt}0\hspace*{6pt}} \\[-1pt]
\hline\\[-13pt]
\multicolumn{1}{|c|}{} &\multicolumn{1}{c}{}&\multicolumn{1}{|@{}c@{}|}{0} \\[-1pt]
\hline\\[-13pt]
\multicolumn{1}{|c|}{0} &\multicolumn{1}{c}{0}&\multicolumn{1}{|@{}c@{}|}{0} \\
\hline
\end{tabular}\hspace*{10pt}
\begin{tabular}{ccc}
\hline\\[-13pt]
\multicolumn{1}{|c|}{0} &\multicolumn{1}{c}{}&\multicolumn{1}{|c@{}|}{\hphantom{0\hspace*{5.8pt}}} \\[-1pt]
\hline\\[-13pt]
\multicolumn{1}{|c|}{} &\multicolumn{1}{c}{}&\multicolumn{1}{|@{}c@{}|}{ } \\[-1pt]
\hline\\[-13pt]
\multicolumn{1}{|c|}{0} &\multicolumn{1}{c}{0}&\multicolumn{1}{|@{}c@{}|}{ } \\
\hline
\end{tabular}\hspace*{10pt}
\begin{tabular}{ccc}
\hline\\[-13pt]
\multicolumn{1}{|c|}{0} &\multicolumn{1}{c}{0}&\multicolumn{1}{|@{}c@{}|}{\hspace*{6pt}0\hspace*{6pt}} \\[-1pt]
\hline\\[-13pt]
\multicolumn{1}{|c|}{} &\multicolumn{1}{c}{0}&\multicolumn{1}{|@{}c@{}|}{0} \\[-1pt]
\hline\\[-13pt]
\multicolumn{1}{|c|}{0} &\multicolumn{1}{c}{0}&\multicolumn{1}{|@{}c@{}|}{0} \\
\hline
\end{tabular}
\end{center}
%
%
\end{example}

\section{Algorithms for extended maximum likelihood estimation}\label
{sec:algo}

In the supplementary material [\citet{FRsupp}], we apply the theory
developed in this article to develop efficient algorithms for extended
maximum likelihood estimation in log-linear models under Poisson and
product multinomial schemes [for which the key integrality assumption
(A2) is satisfied] that are applicable to high-dimensional models and
large tables. Some of these algorithms are implemented in a {\tt
MATLAB} toolbox available at
\url{http://www.stat.cmu.edu/\textasciitilde arinaldo/ExtMLE/}. The final output of
our procedure is the set of estimable mean value and natural
parameters.

\begin{appendix}

\section{Extended exponential families}\label{app:ext.expo}

In this appendix we provide a brief review of the theory of extend
families and its relevance for log-linear models. Along with classic
references on exponential families [\citet{BARN78}, \citet{BRW86},
\citet{CENCOV82}, \citet{LETAC}] and generalizations by \citeauthor{CM01}
(\citeyear{CM01,CM03,CM05,CM08}), we refer the reader to
\citet{ERGM09} and \citet{GEYER09} for treatments more directly
relevant to our problem.

Consider a log-linear model under conditional Poisson sampling scheme
specified by a sampling matrix $\m{V}$ of rank $m$ and a design matrix
$\m{A}$ of the form~\eqref{eq:A}, where $\m{B}$ is of full-rank $d -
m$. Then [see equation~\eqref{eq:qtheta}], the distribution of the
sufficient statistic $\mathbf{z} = \m{B}^\top\mathbf{n}$ form an
exponential family of distributions $\mathcal{E}_{C_{\mathrm{V}}}$ on
$\mathbb{R}^{d-m}$ with densities
\[
q_{\theta}(\mathbf{z}) = \exp \bigl( ( \mathbf{z}, \bftheta) - \psi
(\bftheta)  \bigr),\qquad  \bftheta\in\Theta,
\]
with respect to the base measure $\mu_{\m{V}} = \nu_{\m{V}}^{-1}
\m{B}$, and parameter space $\Theta = \mathbb{R}^{d-m}$. The
\textit{convex support} $C_{\mathrm{V}}$ of $\mathcal{E}_{C_{\mathrm{V}}}$ is the closure
of the convex hull of the support of $\mu_{\m{V}}$. In particular,
$\mathrm{P}$ is a full-dimensional polyhedron in $\mathbb{R}^{d-m}$
and, for every face $F$ of $C_{\mathrm{V}}$, $F$ is the convex hull of some
points in the support of $\mu_{\m{V}}$. Given a realization
$\mathbf{z}$ of the sufficient statistics, the random set\looseness=-1
\begin{equation}\label{eq:mle}
\widehat{\bftheta}(\mathbf{z}) = \widehat{\bftheta} =  \Bigl\{
\bftheta^* \in\Theta\dvtx  q_{\bftheta^*}(\mathbf{z}) = \sup
_{\bftheta\in\Theta} q_{\bftheta}(\mathbf{z})  \Bigr\}
\end{equation}\looseness=0
is the \textit{maximum likelihood estimator} (MLE) of $\bftheta$. If
$\widehat{\bftheta} = \varnothing$ we say that the MLE does not exist.
Existence of the MLE is fully characterized by the geometry of $C_{\mathrm{V}}$,
as the following well-known result indicates; see, for example, Theorem
5.5 in \citet{BRW86} or Theorem 9.13 in \citet{BARN78}.\looseness=-1

\begin{theorem}\label{thm:mle} For a minimal and full
exponential family, the MLE $\widehat{\bftheta}$ exists and is unique
if and only if $\mathbf{z} \in\operatorname{ri}(\mathrm{P})$.
\end{theorem}

Setting $\bolds{\xi}(\bftheta)= \int_{\mathbb{R}^{d-m}} \mathbf{z}
q_{\bftheta}(\mathbf{z})\,d \mu_{\m{V}}(\mathbf{z})$, because of the
minimality of $\mathcal{E}_{C_{\mathrm{V}}} $, one obtains the fundamental
identity $\nabla\psi(\bftheta) = \bolds{\xi}(\bftheta),   \forall
\bftheta\in\Theta,$ where $\nabla$ indicates the gradient.
In particular, if the MLE exists, it satisfies the equation $\widehat
{\bftheta} = (\nabla\psi)^{-1}(\mathbf{z}),$ which is equivalent to the
moment equation $\bolds{\xi}(\widehat
{\bftheta}) = \mathbf{z}$.\vadjust{\goodbreak} 

For any proper face $F$ of $C_{\m{V}}$, let $\mu^F_{\m{V}}$ be the
restriction of $\mu_{\m{V}}$ to $F$. Then, $\mu^F_{\m{V}}$
determines a
new exponential family of distributions, $\mathcal{E}_F$, with
densities with respect to $\mu^F_{\m{V}}$ given by
\[
q^F_{\bftheta}(x) = \exp \bigl( (\mathbf{z,} \bftheta) - \psi
^F(\bftheta)
 \bigr),\qquad  \bftheta\in\Theta_F,
\]
where the natural parameter space is $\Theta_F = \{ \bftheta\in
\Theta\dvtx \exp ( \psi^F(\bftheta)  )< \infty\} = \Theta$, with
$\psi^F(\bftheta) = \log\int_{\mathbb{R}^{d-m}} \exp((\mathbf{z},
\bftheta) ) \,d \mu^F_{\m{V}}(\mathbf{z})$. The convex support of this
new family is $F$, and the existence result of Theorem \ref{thm:mle}
carries over: the MLE exists if and only if the observed sample
$\mathbf{z}$ belongs to $\operatorname{ri}(F)$. However, since~$\mathcal{E}_F$ is supported on a lower-dimensional affine subspace of
$\mathbb{R}^{d-m}$ of dimension $d_F = \operatorname{dim}(F)$, it is no
longer minimal, hence it is unidentifiable.
Nonetheless, if $\mathbf{z} \in\operatorname{ri}(F)$, the MLE of
$\bftheta$ is the set consisting of those $\bftheta$ satisfying the
first order optimality conditions
\begin{equation}\label{eq:mle.F}
\mathbf{z} = \nabla\psi_F(\bftheta)  \quad \mbox{or,
equivalently,}\quad
 \bolds{\xi}^F(\bftheta) = \mathbf{z},
\end{equation}
where $\bolds{\xi}^F(\bftheta)= \int_{\mathbb{R}^{d-m}} \mathbf{z}
q^F_{\bftheta}(\mathbf{z}) \,d \mu^F_{\m{V}}(\mathbf{z})$.


The collection of distributions
\[
\mathcal{E} = \bigcup_F \mathcal{E}_F
\]
as $F$ ranges over all the faces of $C_{\m{V}}$, including $C_{\m{V}}$
itself, is known as the \textit{extended exponential family} of
distributions. With respect to such family $\mathcal{E}$, for any
observed sample $\mathbf{z}$, the MLE, or \textit{extended MLE}, is
always well defined and is the set of solutions to \eqref{eq:mle.F},
where $F$ is the unique face containing $\mathbf{z}$ in its relative
interior.




\section{Proofs}\label{appendix.proof}

This appendix contains the proofs of some results stated in the
article. The remaining proofs can be found in the supplementary
material \citet{FRsupp}. Throughout, we assume familiarity with basic
notions of polyhedral geometry; see \citet{ZIE98}, \citet{SCH98} and
\citet{ROCK70} for in-depth treatments, and Section 2.1 of
\citet{ERGM09} for a brief review of the concepts directly relevant to
our setting.

\begin{pf*}{Proof of Theorem \ref{thm:ex.mle}}
We first assume that $\m{A}$ is of full rank $d$. If $\mathcal{N} = \{
\mathbf{0} \}$, then the convex support is the $d$-dimensional
polyhedral cone~$C_{\m{A}}$, so the result follows directly from
Theorem \ref{thm:mle}. Thus, throughout the remainder of the proof we
consider the case $0 < \operatorname{dim}(\mathcal{N}) < d$. For now,
we further assume that $\m{A}$ is of the form \eqref{eq:A}.

By standard minimality arguments, we can work with the exponential
family supported on $S = \{ \mathbf{z} \dvtx \mathbf{z} =
\m{B}^\top\mathbf{x}, \mathbf{x} \in\mathbb{N}^{\mathcal{I}},
\m{V}^\top\bfx= 1\}$. By assumption~(A1), the convex support
$C_{\m{V}}$, which is the closure of the set $\operatorname{conv}(S)$,
is a~full-dimensional polyhedron in $\mathbb{R}^{d-m}$. In particular,
the parameter space is $\mathbb{R}^{d-m}$. The MLE exists and is
unique\vadjust{\goodbreak}
if and only if $\mathbf{z} \in\operatorname{ri}(C_{\m{V}})$ by Theorem
\ref{thm:mle}. We now show that this happens if and only if $\mathbf{t}
\in\operatorname{ri}(C_\m{A})$. We first use the integrality assumption
(A2) to obtain a simpler representation of $C_{\m{V}}$.

\begin{lemma}\label{lem:CN}
\[
C_{\m{V}}= \{ \m{B}^\top\bfx\dvtx \bfx\in\mathbb{R}^{\mathcal
{I}}_{\geq0}, \m{V}^\top\bfx= \mathbf{1} \}.
\]
\end{lemma}

\begin{pf}
Since $ \{ \m{B}^\top\bfx\dvtx \bfx\in
\mathbb{R}^{\mathcal{I}}_{\geq0}, \m{V}^\top\bfx= 1 \}$ is a polyhedron
(hence closed and convex), it must contain $C_{\m{V}}$. To show the
reverse inclusion, let $\bfz^* \in\{ \m{B}^\top\bfx\dvtx
\bfx\in\mathbb{R}^{\mathcal{I}}_{\geq0}, \m{V}^\top\bfx= \mathbf{1}
\}$. Then, $\bfz^* = \m{B}^\top\bfx^*$ for some $\bfx^* \in\{ \bfx
\dvtx \bfx\in\mathbb{R}^{\mathcal{I}}_{\geq0},\allowbreak \m{V}^\top\bfx=
\mathbf{1} \}$. By the integrality assumption (A2),
\[
\bfx^* \in\operatorname{conv}(\{ \bfx\dvtx \bfx\in\mathbb
{N}^{\mathcal{I}}, \m{V}^\top\bfx= \mathbf{1} \}),
\]
which by linearity implies that $\bfz^*
\in\operatorname{conv}(\m{B}^\top \bfx\dvtx
\bfx\in\mathbb{N}^{\mathcal{I}}, \m{V}^\top\bfx= \mathbf{1} \})
\subseteq C_{\m{V}}$, as claimed.
\end{pf}

For design matrices of the form \eqref{eq:A}, the claim in the theorem
follows directly from the next lemma.

\begin{lemma}\label{ref:homo}
There exists a homomorphism from the face lattice of $C_{\m{V}}$ to the
face lattice of $C_{\m{A}}$ that associates to each face of $C_{\m{V}}$
of dimension $d_F$ the (unique) face of $C_{\m{A}}$ of dimension $m +
d_F$ containing it.
\end{lemma}

\begin{pf}
Instead of concerning ourselves with $C_{\m{V}}$, we find it convenient
to deal with the $d - m$-dimensional polyhedron in $\mathbb{R}^d$
\begin{equation}\label{eq:TN}\quad
T_{\m{V}} = C_{\m{A}} \cap\{ \bft= (t_1,\ldots,t_d)^\top\in
\mathbb{R}^d \dvtx  t_j = 1, j = d-m + 1,\ldots, d\}.
\end{equation}
%
In light of the next result, $T_{\m{V}}$ and $C_{\m{V}}$ have the same
combinatorial properties.

\begin{lemma}\label{lem:CN.TN}
The polyhedra $T_{\m{V}}$ and $C_{\m{V}}$ are combinatorially
equivalent. 
\end{lemma}

\begin{pf}
By Lemma \ref{lem:CN}, $\bfz\in C_{\m{V}}$ if and only if
$
 (
{{\bfz} \atop {\mathbf{1}}}
 ) \in T_{\m{V}}.
$ Thus the coordinate projection map $\pi\dvtx \mathbb{R}^d
\rightarrow\mathbb{R}^{d-m}$ given by $\pi(x_1,\ldots,x_d) =
(x_{1},\ldots,x_{d-m})$ defines a bijection between $C_{\m{V}}$ and
$T_{\m{V}}$. Since $\pi$ is a linear mapping, $C_{\m{V}}$ and~$T_{\m{V}}$ are affinely equivalent, hence combinatorially equivalent.
\end{pf}

It follows from Lemma \ref{lem:CN.TN} that there exists a bijection
between $C_{\m{V}}$ and~$T_{\m{V}}$ that is also a bijection between
boundary points of $C_{\m{V}}$ and points on the relative boundary of~$T_{\m{V}}$ in such a way that the face lattices of $C_{\m{V}}$ and
$T_{\m{V}}$ are identical. Note also that isomorphic faces of the
polyhedra have the same dimension. Therefore, it is sufficient to prove
that the claim of the theorem holds for $T_{\m{V}}$ instead of
$C_{\m{V}}$.

Using the $\mathcal{H}$-representation [see, e.g., \citet{ZIE98},
\citet{SCH98}] we write
\begin{equation}\label{eq:H.CA}
C_{\m{A}} = \{ \bft\in\mathbb{R}^d \dvtx \m{C} \bft\leq\mathbf{0}\}\vadjust{\goodbreak}
\end{equation}
for some matrix $\m{C}$, where we can assume that no inequality is
redundant. In particular, any face $F$ of $C_{\m{A}}$ of co-dimension
$k$ can be written as
\[
\{ \mathbf{t} \dvtx \m{C} \mathbf{t} \leq\mathbf{0},
(\mathbf{c}_j,\bft) = 0,
j=1,\ldots,k\},
\]
where $(\mathbf{c}_1,\ldots,\mathbf{c}_k)$ are the $k$ rows of $\m{C}$
that define the $k$ supporting hyperplanes whose intersection with
$C_{\m{A}}$ is precisely $F$. Define
\[
\m{T} =  \left[\matrix{ \m{0}& \m{I}_m  }\right],
\]
where $\m{0}$ is the $ m \times(d-m)$ matrix of zeros, and $\m{I}_m$ is
the $m \times m$ identity matrix. Thus, $T_{\m{V}}$ is the set of
points in $\mathbb{R}^d$ given by $\{ \mathbf{t} \dvtx \m{D} \bft\leq
\mathbf{b}\}$, with
\[
\m{D} =  \left[
\matrix{ \hspace*{-2pt}\phantom{-}\m{C}'\cr
\phantom{-}\m{T}\cr
-\m{T}
}
 \right]  \quad\mbox{and} \quad  \mathbf{b} =  \left[
\matrix{
\phantom{-}\mathbf{0}\cr
\phantom{-}\mathbf{1}\cr
-\mathbf{1}
}
 \right],
\]
where $\m{C}'$ is the sub-matrix of $\m{C}$ obtained by removing the
rows corresponding to inequalities that may have become redundant once
the sampling constraint are enforced. These inequalities are the
precisely the defining inequalities for the facets that do not
intersect the affine space $\{ \bft\dvtx \m{T} \bft= 1\}$. Notice
that, by (A1), the dimension of $T_{\m{V}}$ is equal to $d$ minus the
rank of
\[
 \left[
\matrix{
\phantom{-}\m{T}\cr
-\m{T}
}
 \right],
\]
which is $m$. Next, any face $F$ of $T_{\m{V}}$ of co-dimension $k$ can
be written as
\[
F = \{ \bft\dvtx \m{D} \bft\leq\mathbf{b}, (\mathbf{d}_j,\bft) = 0,
j=1,\ldots,k \},
\]
where $(\mathbf{d}_1,\ldots,\mathbf{d}_k)$ are the $k$ rows of $\m{C}'$
that define the $k$ supporting hyperplanes of $F$. Since the points in
$F$ satisfy all the inequalities \eqref{eq:H.CA}, it follows that $F$
is contained in the set $F' = \{ \bft\dvtx \m{C} \bft\leq\mathbf{0},
(\mathbf{d}_j,\bft) = 0, j=1,\ldots,k\}$, which is a face of
$C_{\m{A}}$ of co-dimension $k$. It is also immediate to see that~$F'$
is the smallest such face. Furthermore, if $G$ is a different face of
$T_{\m{V}}$ of co-dimension~$k$, it is defined by a different set of
equalities, so it is contained in a~different face of $C_{\m{A}}$ (of
co-dimension $k$). If $G$ is instead of co-dimension $k' > k$ and is
also a face of $F$, then, $G = \{ \bft\dvtx \m{D} \bft\leq\mathbf{b},
(\mathbf{d}_j,\bft) = 0$, $j=1,\ldots,k,\ldots,k' \}$, so that $G$ is
contained in the set $\{ \bft\dvtx \m{C} \bft\leq\mathbf{0},
(\mathbf{d}_j,\bft) = 0$, $j=1,\ldots,k'\}$, which is a face of
$C_{\m{A}}$ of co-dimension $k'$ and also a face of~$F'$.

Therefore, the mapping that associates to each face of $T_{\m{V}}$ the
smallest face of $C_{\m{A}}$ containing it (and of the same
co-dimension) is a lattice homomorphism from the face lattice of
$T_{\m{V}}$ to the face lattice of $C_{\m{A}}$. Furthermore, since the
homomorphism just described is between faces of the same co-dimension,
and $\operatorname{dim}(T_{\m{V}}) = d - m$ while
$\operatorname{dim}(C_{\m {A}}) = d$, each face of $T_{\m{V}}$ of
dimension $d_F$ is mapped to a face of $C_{\m{A}}$ of dimension $m +
d_F$.
\end{pf}

Thus far we have assumed that the design matrix $\m{A}$ is of full rank
and has the form specified by equation \eqref{eq:A}. Now let $\m{A}'$
be any design matrix with row span $\mathcal{M}$, not necessarily of
the form \eqref{eq:A}, or not even of full rank. Then, $C_{\m{A}'}$ is
also a polyhedral cone of dimension $d$, though its ambient dimension
may be larger. As $\m{A}'$ and $\m{A}$ have the same null space, the
cones~$C_{\m{A}'}$ and $C_{\m{A}}$ are affinely isomorphic, hence
combinatorially equivalent. Thus, $\mathbf{t}' = (\m{A}')^\top x \in
\operatorname{ri}(C_{\m{A}'})$ if and only if $\mathbf{t} =
\m{A}^\top\mathbf{x}$, which shows that the theorem holds for any
generic design matrix $\m{A}$.
\end{pf*}

\begin{pf*}{Proof of Theorem \ref{thm:main}}
We show that, under both Poisson and product multinomial scheme, the
MLE exists, is unique and is identical in both cases if and only if
$\mathbf{t} = \m{A}\mathbf{n}$ is a point in the relative interior of
$C_{\m{A}}$. If $\mathbf{t}$ belongs to the relative interior of a face
$F$, then both log-likelihood functions realize their suprema along
sequences of points $\bolds{\mu}_n \subset\mathcal{M}$ for which the
limit $\exp(\bolds{\mu}_n) = \widehat{\mathbf{m}}{}^{\m{e}}$ is unique,
satisfies the moment equations $\Pi_{\mathcal{M}}\mathbf{n} =
\Pi_{\mathcal{M}} \widehat{\mathbf{n}}$ and
$\operatorname{supp}(\widehat{\mathbf{m}}) = \mathcal{F}$.

First, we consider the problem of maximizing the log-likelihood
$\ell^P(\bfmu) = (\mathbf{n}, \bfmu) - \sum_{i \in\mathcal{I}}
\exp(\bfmu(i))$ under Poisson sampling scheme. Suppose $\mathbf{t} =
\m{A}^\top\mathbf{n}$ lies inside the relative interior of a proper
face $F$ of $C_{\m{A}}$ with corresponding facial set~$\mathcal{F}$.
Then, there exists a $\mathbf{z}_{F} \in\operatorname{kernel}(\m{A}) =
\mathcal{M}^\bot$ such that the vector $\mathbf{x}_{F} =
\mathbf{n}+\mathbf{z}_{F}$ satisfies $\mathbf{t} =
\m{A}^\top\mathbf{x}_F$ and
$\operatorname{supp}(\mathbf{n}+\mathbf{z}_{F}) = \mathcal{F}$.
Furthermore, since, for any $\bfmu\in\mathcal{M}$,
$(\mathbf{z}_{F},\bfmu)=0$, $\ell^P(\bfmu) = (\mathbf{x}_{F}, \bfmu) -
\sum_{i \in\mathcal{I}} \exp(\bfmu(i))$.

Define $\ell^P_{\mathcal{F}}$ and $\ell^P_{\mathcal{F}^c}$ to be the
restriction of $\ell^P$ on $\pi_{\mathcal{F}}(\mathcal{M})$ and
$\pi_{\mathcal{F}^c}(\mathcal{M})$, respectively. Explicitly,
\[
\ell^P_{\mathcal{F}}(\bolds{\mu}) = (\mathbf{x}_F, \pi_{\mathcal
{F}}(\bfmu)) - \sum_{i \in\mathcal{F}} \exp(\bfmu(i)) =
(\mathbf{x}_{F}, \bfmu) - \sum_{i \in\mathcal{F}} \exp(\bfmu(i))
\]
and $\ell^P_{\mathcal{F}^c}(\bolds{\mu}) = -\sum_{i \in\mathcal{F}^c}
\exp(\bfmu(i) )$. Therefore, $\ell^P(\bolds{\mu})
=\ell^P_{\mathcal{F}}(\bolds{\mu}) + \ell^P_{\mathcal{F}^c}(\bolds
{\mu})$. On $\pi_{\mathcal{F}}(\mathcal{M})$, the function
$\ell^P_{\mathcal{F}}$ is bounded from above, continuous and strictly
concave, so it is maximized by the unique point $\bfmu_{\mathcal{F}}^*
\in\pi_{\mathcal{F}}(\mathcal{M})$ that satisfy the first order
optimality conditions on the differential of $\ell^P_{\mathcal{F}}$
[see \citet{HAB74}, Chapter 2] given by
\begin{equation}\label{eq:diff}\qquad
 ( \bflambda_{\mathcal{F}}, \exp(\bfmu_{\mathcal{F}}^*)
 ) =  ( \bflambda_{\mathcal{F}}, \pi_{\mathcal{F}}(\bfx _F)
 ) =  ( \bflambda_{\mathcal{F}}, \bfn )\qquad
\forall\bflambda_{\mathcal{F}} \in\pi_{\mathcal{F}}(\mathcal{M}),
\end{equation}
where the second equality holds since $\mathbf{x}_F
\in\mathcal{M}^\bot$ and $\operatorname{supp}(X_\mathcal{F}) =
\mathcal{F}$.

On the other hand, on $\pi_{\mathcal{F}^c}(\mathcal{M})$, the function
$\ell^P_{\mathcal{F}^c}$ is negative and strictly decreasing in each
coordinate of its argument. Thus,
\[
\sup_{\bfmu\in\mathcal{M} } \ell^P_{\mathcal{P}(\bfmu)} \leq
\sup_{\bfmu_{\mathcal{F}} \in\pi_{\mathcal{F}}(\mathcal{M})}
\ell^P_{\mathcal{F}}(\bfmu_{\mathcal{F}}) = \ell^P_{\mathcal
{F}}(\bfmu^*_{\mathcal{F}}).
\]

We now show that the above inequality is in fact an equality by
finding a sequence $\{ \bfmu_n \} \subset\mathcal{M}$ such that
\[
\lim_n \ell^P(\bfmu_n) = \ell^P_{\mathcal{F}}(\bfmu^*_{\mathcal{F}}).
\]
To this end, let $\bfmu^*$ be any vector in $\mathcal{M}$ such that
$\pi_{\mathcal{F}}(\bfmu^*) = \bfmu^*_{\mathcal{F}}$. Next, since~$\mathcal{F}$ is a facial set, there exists a sequence $\{ \bfgamma_n
\} \subset\mathcal{M}$ such that:
\begin{enumerate}[(ii)]
\item[(i)] if $i \in\mathcal{F}$, then $\bfgamma_n(i) = 0$, for all
$n$;
\item[(ii)] if $i \in\mathcal{F}^c$, then $\bfgamma_n(i) < 0$
for all $n$ and $\lim_n \bfgamma_n(i) = - \infty$ (the rate at which
these series diverge to infinity being arbitrarily fast).
\end{enumerate}
Define the sequence $\{ \bfmu_n \} \subset\mathcal{M}$ as $\bfmu_n
= \bfmu^* + \bfgamma_n$. Then,
\[
\lim_n \bfmu_n(i) =  \cases{
\bfmu^*(i) &\quad if $i \in\mathcal{F}$,\cr
-\infty&\quad if $i \in\mathcal{F}^c$,
}
\]
from which it follows that
\begin{eqnarray*}
\lim_n \ell^P(\bfmu_n) &=& \lim_n \ell^P_{\mathcal{F}}(\pi
_{\mathcal{F}}(\bfmu_n)) + \lim_n \ell^P_{\mathcal{F}^c}(\pi
_{\mathcal{F}^c}(\bfmu_n)) \\
&=& \ell^P_{\mathcal{F}}(\bfmu
^*_{\mathcal{F}}) + \lim_n \ell^P_{\mathcal{F}^c}(\pi_{\mathcal
{F}^c}(\bfmu_n)) = \ell^P_{\mathcal{F}}(\bfmu^*_{\mathcal{F}}),
\end{eqnarray*}
as desired, since
\[
\lim_n \ell^P_{\mathcal{F}^c}(\pi_{\mathcal{F}^c}(\bfmu_n)) =
\sum_{i \in\mathcal{F}^c} \lim_n \exp(\bfmu_n(i)) = 0.
\]
Set $\widehat{\bfm}{}^{\m{e}} = \lim_n \exp(\bfmu_n)$, and notice that
$\widehat{\bfm}{}^{\m{e}}$ is the unique vector in
$\mathbb{R}^{\mathcal{I}}$ such that
\[
 \cases{
\pi_{\mathcal{F}}(\widehat{\bfm}{}^{\m{e}})  =  \exp(\bfmu
^*_{\mathcal{F}}),\cr
\pi_{\mathcal{F}^c}(\widehat{\bfm}{}^{\m{e}})  =  \mathbf{0},
}
\]
where uniqueness stems from the uniqueness of $\bfmu^*_{\mathcal{F}}$
(it is clear that, while~$\widehat{\bfm}{}^{\m{e}}$ is unique, the
sequence $\{ \bfmu_n \}$ is not). Furthermore,
$\widehat{\bfm}{}^{\m{e}}$ is random, as it depends on the facial set
$\mathcal{F}$ associated to the face of $C_{\m{A}}$ exposed by
$\mathbf{t} = \m{A}^\top\bfn$. Finally, in virtue of the fact that
$\operatorname{supp}(\bfn) \subseteq \mathcal{F}$, we see that, for any
$\bflambda\in\mathcal{M}$,
\[
 ( \bflambda, \widehat{\bfm}{}^{\m{e}}  ) =  (
\bflambda_{\mathcal{F}}, \exp(\mu_{\mathcal{F}}^*)  ) \quad\mbox{and}\quad
   ( \bflambda_{\mathcal{F}}, \bfn ) =  ( \bflambda,
\bfn )
\]
so that, using \eqref{eq:diff}, $\widehat{\bfm}{}^{\m{e}}$ can be
characterized as the unique point in $\overline{M}$ such that
\[
 ( \bflambda, \widehat{\bfm}{}^{\m{e}}  ) =  (
\bflambda, \bfn ) \qquad \forall\bflambda\in\mathcal{M},
\]
or, equivalently,
\begin{equation}\label{eq:unique}
\m{A}^\top\widehat{\bfm}{}^{\m{e}} = \m{A}^\top\bfn \quad\mbox {or}\quad
 \Pi_{\mathcal{M}} \widehat{\bfm}{}^{\m{e}} = \Pi _{\mathcal{M}}
\bfn.
\end{equation}

If we instead want to maximize the log-likelihood function $\ell^M$
under product multinomial sampling, we need to consider only the points
$\tilde{\bolds{\mu}}$ inside $\widetilde{\mathcal{M}}$ as in equation
(\ref{eq:Mtilde}). Fortunately, this restriction is inconsequential.
First note that, by \eqref{eq:unique} and because $\mathcal{N} \subset
\mathcal{M}$, the limit $\bolds{\mu}^*$ satisfies the constraints $\{
(\bolds{\chi}_j, \exp(\bolds{\mu}^*)) = N_j, j=1,\ldots,r\}$.\vadjust{\goodbreak} Next,
since $\ell^M$ and $\ell^P$ differ by a constant on
$\widetilde{\mathcal{M}}$ and $\widetilde{\mathcal{M}}
\subset\mathcal{M}$, we have that
\[
\ell^M(\bolds{\mu}^*) = \sup_{\tilde{\bolds{\mu}} \in\widetilde
{\mathcal{M}}} \ell^M(\tilde{\bolds{\mu}}).
\]
We conclude that the log-likelihood functions under both the Poisson
and product multinomial model must have the same maximizer
$\widehat{\mathbf{m}}$.

Finally, we note that if $\mathbf{t} \in\operatorname{ri}(C_{\m{A}})$,
so that $\mathcal{F} = \mathcal{I}$, the arguments simplify.
Explicitly, there exists a point $\bolds{\mu}^*
\in\widetilde{\mathcal{M}} \subset \mathcal{M}$ such that
\begin{eqnarray*}
\sup_{\bolds{\mu} \in\mathcal{M}} \ell^P (\bolds{\mu}) & = & \ell
^P (\bolds{\mu}^*),\\
\sup_{\tilde{\bolds{\mu}} \in\widetilde{\mathcal{M}}} \ell^M
(\tilde{\bolds{\mu}}) & = & \ell^M (\bolds{\mu}^*)
\end{eqnarray*}
which we can obtain as the unique point $\widehat{\bfm}{}^{\m{e}} \in
\overline{M}$ with $\operatorname{supp}(\widehat{\bfm}{}^{\m{e}}) =
\mathcal{I}$ satisfying~\eqref{eq:unique}.
\end{pf*}
\end{appendix}

\section*{Acknowledgment}
The authors would like to thank the Associate Editor for invaluable help and support.


%

\printaddresses

\end{document}